\RequirePackage{ifpdf}
\ifpdf 
\documentclass[pdftex]{sigma}
\else
\documentclass{sigma}
\fi

\usepackage{euscript}

\DeclareMathOperator{\Imp}{\mathrm{Im}\,}
\DeclareMathOperator{\Rep}{\mathrm{Rep}\,}
\DeclareMathOperator{\End}{\mathrm{End}}
\DeclareMathOperator{\Hom}{\mathrm{Hom}}
\DeclareMathOperator{\Idem}{\mathrm{Idem}}
\DeclareMathOperator{\Diag}{\mathrm{diag}}
\DeclareMathOperator{\Mor}{\mathrm{Mor}}

\def\Id{\mathrm{Id}\,}

\begin{document}
\allowdisplaybreaks

\renewcommand{\PaperNumber}{042}

\FirstPageHeading

\ShortArticleName{On Transitive Systems of Subspaces in a Hilbert
Space}

\ArticleName{On Transitive Systems of Subspaces in a Hilbert
Space}

\Author{Yuliya P. MOSKALEVA~$^\dag$ and Yurii S. SAMO\v{I}LENKO~$^\ddag$}

\AuthorNameForHeading{Yu.P. Moskaleva and Yu.S. Samo\v\i{}lenko}

\Address{$^\dag$ Taurida National University, 4 Vernads'kyi Str.,
  Simferopol, 95007 Ukraine}

\EmailD{\href{mailto:YulMosk@mail.ru}{YulMosk@mail.ru}}

\Address{$^\ddag$~Institute of Mathematics, National Academy of  Sciences of Ukraine,\\
$\phantom{^\ddag}$~3 Tereshchenkivs'ka Str., Kyiv-4, 01601 Ukraine}

\EmailD{\href{mailto:yurii_sam@imath.kiev.ua}{yurii\_sam@imath.kiev.ua}}

\ArticleDates{Received February 27, 2006; Published online April 12, 2006}

\bigskip

\leftline{\small \bfseries \itshape Submitted by Anatoly Klimyk}

\Abstract{Methods of $*$-representations in Hilbert space are
applied to study of systems of $n$ subspaces in a linear space. It
is proved that the problem of description of
 $n$-transitive subspaces in a f\/inite-dimensional linear space is $*$-wild for $n \geq 5$.}


\Classification{47A62; 16G20}

\section{Introduction}

Systems of $n$ subspaces $H_1,H_2,\ldots,H_n$ of a Hilbert space
$H$, denoted in the sequel by $S=(H;H_1,H_2,\ldots,H_n)$, is a
mathematical object that traditionally draws an interest both by
itself~\cite{Ha1,GP,N,EW} and in connection with the discussion on
whether there exists a deeper connection between this object and
the famous H.~Weyl problem, the Coxeter groups, singularity
theory, and physical applications.

Systems of subspaces that can be regarded as candidates for being
the simplest building blocks for arbitrary systems of subspaces
are those that are indecomposable or transitive~\cite{GP,N,EW}.
A~description of transitive and indecomposable systems is carried
out up to an isomorphism of the systems of subspaces. For a
description of transitive and indecomposable systems of two
subspaces of a Hilbert space, as well as for transitive and
indecomposable triples of a f\/inite dimensional linear space,
see, e.g.,~\cite{EW}. For an inf\/inite dimensional space, not
only the problem of description but even the problem of existence
of transitive and indecomposable triples of subspaces is an
unsolved problem~\cite{Ha2}.  For a f\/inite dimensional linear
space, transitive quadruples of subspaces are described
in~\cite{B}, and \cite{GP,N} give indecomposable quadruples.
Examples of nonisomorphic transitive and indecomposable systems of
four subspaces in an inf\/inite dimensional space can be found,
e.g., in~\cite{EW}.

In~\cite{EW} the authors make a conjecture that there is a
connection between systems of $n$ subspaces and representations of
$*$-algebras that are generated by the projections, --- ``There
seems to be interesting relations of systems of $n$-subspaces with
the study of representations of $*$-algebras generated by
idempotents by S. Kruglyak, V. Ostrovskyi, V.  Rabanovich, Yu.
Samo\v\i{}lenko and other. But we do not know the exact
implication \dots''. The present article deals with this
implication.

Let us consider systems of subspaces of the form
$S_\pi=(H;P_1H,P_2H,\ldots,P_nH)$, where the orthogonal
projections $P_1, P_2, \ldots, P_n$ make a $*$-representation
$\pi$ of the $*$-algebra genera\-ted by the projections, and $H$
is the representation space. For the \hbox{$*$-algeb}ras
$\EuScript P_{4,{\rm com}}=\mathbb C\langle p_1,p_2,p_3, p_4 \,
|\, p^2_k=p^*_k=p_k, \big
[\sum\limits_{k=1}^4p_k,p_i\big]=0,\,\forall\, i=1,2,3,4\rangle$,
it was proved in~\cite{MS} that irreducible inequivalent
$*$-rep\-resentations $\pi$ of the $*$-algebra $\EuScript
P_{4,{\rm com}}$ make a complete list of nonisomorphic transitive
quadruples of subspaces $S_\pi$ of a f\/inite dimensional linear
space.

In this paper, we make an analysis of complexity of the
description problem for transitive systems of subspaces
$S=(H;H_1,H_2,\ldots,H_n)$ for $n\geq5$. In Section~3, we prove
that it is an extremely dif\/f\/icult problem to describe
nonisomorphic transitive quintuples of subspaces
$S=(H;P_1H,P_2H,\ldots,P_5H)$ even under the assumption that the
sum of the corresponding f\/ive projections equals $2I$; in other
words, the problem of describing inequivalent $*$-representations
of the $*$-algebras that give rise to nonisomorphic transitive
systems, is $*$-wild.

Since the problem of describing the system of $n$ subspaces up to
an isomorphism is complicated, it seems natural to describe
transitive systems that correspond to $*$-representations of
various algebras generated by projections~(Sections~4 and~5).

In Section~4, we consider transitive systems $S_{\pi}$ of $n$
subspaces, where $\pi\in\Rep\mathcal P_{n,\alpha}$, $\mathcal
P_{n,\alpha}=\mathbb C\langle
p_1,p_2,\ldots,p_n\,|\,p_1+p_2+\cdots+p_n=\alpha e, p_j^2=p_j,
p_j^*=p_j, \, \forall \,j=1,\ldots,n\rangle$, and $\alpha$ takes
values in a f\/ixed set. In Section~5, using nonisomorphic
transitive systems $S_{\pi}$ of $n$ subspaces, where $\pi$ belongs
to $\Rep\mathcal P_{n,\alpha}$, we construct nonisomorphic
transitive systems $S_{\hat{\pi}}$ of $n+1$ subspaces, where
$\hat{\pi}$ is in $\Rep\mathcal P_{n,{\rm abo},\tau}$, $\mathcal
P_{n,{\rm abo},\tau}=\mathbb C\langle
q_1,q_2,\ldots,q_n,p\,|\,q_1+q_2+\cdots+q_n=e, q_jpq_j=\tau q_j$,
$q_j^2=q_j,\, q_j^*=q_j,  \forall\, j=1,\ldots,n, p^2=p,
p^*=p\rangle$.

\section[Definitions and main properties]{Def\/initions and main properties}

In this section we make necessary def\/initions and recall known
facts; the proofs can be found in~\cite{EW,OS}. Let $H$ be a
Hilbert space and $H_1, H_2, \ldots, H_n$ be $n$ subspaces of $H$.
Denote by $S=(H;H_1,H_2,\ldots,H_n)$ the system of $n$ subspaces
of the space $H$. Let $S=(H;H_1,H_2,\ldots,H_n)$ be a system of
$n$ subspaces of a Hilbert space $H$ and
$\tilde{S}=(\tilde{H};\tilde{H}_1,\tilde{H}_2,\ldots,\tilde{H}_n)$
a system of $n$ subspaces of a Hilbert space $\tilde{H}$. A linear
map $R:H\rightarrow\tilde{H}$ from the space $H$ to the space
$\tilde{H}$ is called a homomorphism of the system $S$ into the
system $\tilde{S}$ and denoted by $R:S\rightarrow\tilde{S}$, if $
R(H_i)\subset\tilde{H}_i$, $ i=1,\ldots,n $. A homomorphism
$R:S\rightarrow\tilde{S}$ of a system $S$ into a system
$\tilde{S}$ is called an isomorphism, $R:S\rightarrow\tilde{S}$,
if the mapping $R:H\rightarrow\tilde{H}$ is a bijection and
$R(H_i)=\tilde{H}_i$, $\forall\, i=1,\ldots,n $. Systems $S$ and
$\tilde{S}$ will be called isomorphic, denoted by
$S\cong\tilde{S}$, if there exists an isomorphism
$R:S\rightarrow\tilde{S}$.

Denote by $\Hom(S,\tilde{S})$ the set of homomorphisms of a system
$S$ into a system $\tilde{S}$ and by $\End(S):=\Hom(S,S)$ the
algebra of endomorphisms of $S$ into $S$, that is,
\[
\End(S)=\{R\in B(H)\,|\,R(H_i)\subset H_i, i=1,\ldots,n\}.
\]

A system $S=(H;H_1,H_2,\ldots,H_n)$ of $n$ subspaces of a Hilbert
space $H$ is called transitive, if $\End(S)=\mathbb C I_H$.

Denote
\[
\Idem(S)=\{R\in B(H)\,|\,R(H_i)\subset H_i, i=1,\ldots,n,R^2=R\}.
\]

A system $S=(H;H_1,H_2,\ldots,H_n)$ of $n$ subspaces of a space
$H$ is called indecomposable, if $\Idem(S)=\{0,I_H\}$.

Isomorphic systems are either simultaneously transitive or
intransitive, decomposable or indecomposable. We say that
$S\cong\tilde{S}$ up to permutation of subspaces, if there exists
a permutation $\sigma\in S_n$ such that the systems $\sigma(S)$
and $\tilde{S}$ are isomorphic, where
$\sigma(S)=(H;H_{\sigma(1)},H_{\sigma(2)},\ldots,H_{\sigma(n)})$,
so that there exists an invertible operator
$R:H\rightarrow\tilde{H}$ such that
$R(H_{\sigma(i)})=\tilde{H}_i$, $\forall\, i=1,\ldots,n$.

Let us now recall the notion of unitary equivalence for systems
and collections of orthogonal projections. Systems $S$ and
$\tilde{S}$ are called unitary equivalent, or simply equivalent,
if $S\cong\tilde{S}$ and it is possible to choose the isomorphism
$R:S\rightarrow\tilde{S}$ to be a unitary operator.

To every system $S=(H;H_1,H_2,\ldots,H_n)$ of $n$ subspaces of a
Hilbert space $H$, one can naturally associate a system of
orthogonal projections $P_1, P_2, \ldots, P_n$, where $P_i$ is the
orthogonal projection operator onto the space $H_i$,
$i=1,\ldots,n$. A system of projections $P_1, P_2, \ldots, P_n$ on
a Hilbert space $H$ such that $\Imp P_i=H_i$ for $i=1,\ldots,n$ is
called a system of orthogonal projections associated to the system
of subspaces, $S=(H;H_1,H_2,\ldots,H_n)$. Conversely, to each
system of projections there naturally corresponds a system of
subspaces. A system $S=(H;P_1H, P_2H,\ldots, P_nH)$ is called a
system corresponding to the system of projections $P_1,P_2,
\ldots, P_n$.

A system of orthogonal projections $P_1, P_2, \ldots, P_n$ on a
Hilbert space $H$ is called unitary equivalent to a system
$\tilde{P}_1, \tilde{P}_2, \ldots, \tilde{P}_n$ on a Hilbert space
$\tilde{H}$, if there exists a unitary operator
$R:H\rightarrow\tilde{H}$ such that $RP_i=\tilde{P}_iR$,
$i=1,\ldots,n$. Systems $S$ and $\tilde{S}$ are unitary equivalent
if and only if the corresponding systems of orthogonal projections
are unitary equivalent.

A system of orthogonal projections $P_1, P_2, \ldots, P_n$ on a
Hilbert space $H$ is called irreducible if zero and $H$ are the
only invariant subspaces. Unitary equivalent systems of orthogonal
projections are both either reducible or irreducible.

If systems $S$ and $\tilde{S}$ are unitary equivalent, then
$S\cong\tilde{S}$. The converse is not true.

\begin{example}
  Let $S=(\mathbb C^2;\mathbb C(1,0),\mathbb
  C(\cos\theta,\sin\theta))$, $\theta\in(0,\pi/2)$, and
  $\tilde{S}=(\mathbb C^2;\mathbb C(1,0),\mathbb C(0,1))$. The
  decomposable system $S$ that  corresponds to an irreducible
  pair of orthogonal projections, is isomorphic but not unitary
  equivalent to the decomposable system $\tilde{S}$ that corresponds
  to a reducible pair of orthogonal projections.
\end{example}

Finally, let us mention the relationship between the notions of
transitivity, indecomposability, and irreducibility. If a system
of subspaces is transitive, then it is indecomposable, but not
vice versa.  Indecomposability of a system of subspaces implies
irreducibility of the corresponding system of orthogonal
projections, but not conversely.

\section[On $*$-wildness of the description problem for transitive
  systems of $n$ subspaces for $n\geq5$]{On $\boldsymbol{*}$-wildness of the description problem\\ for transitive
  systems of $\boldsymbol{n}$ subspaces for $\boldsymbol{n\geq5}$}

\subsection[On $*$-wildness of the description problem for transitive
  systems that correspond to orthogonal projections]{On $\boldsymbol{*}$-wildness
  of the description problem for transitive systems\\ that correspond to orthogonal projections}

A description of transitive quadruples of subspaces of a f\/inite
dimensional linear space is given in~\cite{B}. We will show that
such a problem for $n$ subspaces, $n\geq5$, is extremely
complicated ($*$-wild).

Consider a system of f\/ive subspaces, which corresponds to the
 f\/ive orthogonal projections
\begin{gather*}
P_1=\begin{pmatrix}I&0\\0&0\end{pmatrix},\quad
P_2=\begin{pmatrix}0&0\\0&I\end{pmatrix},\quad
P_3=\frac12\begin{pmatrix}I&I\\I&I\end{pmatrix},\\
P_4=\frac12\begin{pmatrix}I&U\\U^*&I\end{pmatrix},\quad
P_5=\frac12\begin{pmatrix}I&V\\V^*&I\end{pmatrix}
\end{gather*}
that act on the space $\mathcal H=H\oplus H$, where $H$ is a
Hilbert space and $U$ and $V$ are unitary operators. Denote this
system of subspaces by $S_{U,V}$. So, $S_{U,V}=(\mathcal
H;P_1\mathcal H,P_2\mathcal H,P_3\mathcal H,P_4\mathcal
H,P_5\mathcal H)$. Consider the system
$S_{\tilde{U},\tilde{V}}=(\tilde{\mathcal
  H};\tilde{P}_1\tilde{\mathcal H},\tilde{P}_2\tilde{\mathcal
  H},\tilde{P}_3\tilde{\mathcal H},\tilde{P}_4\tilde{\mathcal
  H},\tilde{P}_5\tilde{\mathcal H})$ that corresponds to the
collection of orthogonal projections $\tilde{P}_1$, $\tilde{P}_2$,
$\tilde{P}_3$, $\tilde{P}_4$, $\tilde{P}_5$ that have the above
type and act on the space $\tilde{\mathcal H}=\tilde{H}\oplus
\tilde{H}$; here $\tilde{H}$ is a Hilbert space and $\tilde{U}$,
$\tilde{V}$ is a pair of unitary operators.

\begin{theorem}\label{t1}
  The system $S_{U,V}$ is transitive if and only if the unitary
  operators $U$, $V$ are irreducible. Also, $S_{U,V}\cong
  S_{\tilde{U},\tilde{V}}$ if and only if the pair of unitary
  operators $U$, $V$ is unitary equivalent to the pair of unitary
  operators $\tilde{U}$, $\tilde{V}$.
\end{theorem}

\begin{proof}
  Denote $H_i= P_i\mathcal H$, $i=1,\ldots,5$. For $H_1$ and
  $H_2$, we have
  \[
  H_1=H\oplus 0,\qquad  H_2=0\oplus H.
  \]
  For $H_3$, $H_4$, and $H_5$, respectively,
  \[
  H_3=\{(x,x)\,|\,x\in H\},\quad H_4=\{(Ux,x)\,|\,x\in H\},\qquad
  H_5=\{(Vx,x)|x\in H\}.
  \]

  Let us prove an auxiliary identity
  \begin{gather}
    \{\mathcal R\in B(\mathcal H,\tilde{\mathcal H})\,|\,\mathcal
    R(H_i)\subset\tilde{H}_i, i=1,\ldots,5\}\nonumber
    \\
    \qquad{} =\{R\oplus R\in
    B(\mathcal H,\tilde{\mathcal H})\, |\, R\in
    B(H,\tilde{H}),RU=\tilde{U}R,RV=\tilde{V}R\}.\label{x1}
  \end{gather}
  The f\/irst three inclusions, $\mathcal R(H_i)\subset\tilde{H}_i$,
  $i=1,2,3$, imply that any operator $\mathcal R$ in
  $B(\mathcal H,\tilde{\mathcal H})$ can be represented as $\mathcal
  R=R\oplus R$, where $R\in B(H,\tilde{H})$. The fourth inclusion,
  $\mathcal R(H_4)\subset\tilde{H}_4$, implies $RU=\tilde{U}R$, and
  the f\/ifth one, $\mathcal R(H_5)\subset\tilde{H}_5$, gives
  $RV=\tilde{V}R$.  The converse implications f\/inish the proof
  of~(\ref{x1}).

  It directly follows from~(\ref{x1}) that $S_{U,V}\cong
  S_{\tilde{U},\tilde{V}}$ if and only if the pair of unitary
  operators $U$, $V$ is similar to the pair of unitary operators
  $\tilde{U}$, $\tilde{V}$. By~\cite{OS}, a pair of unitary operators
  $U$, $V$ is similar to a pair of unitary operators $\tilde{U}$,
  $\tilde{V}$ if and only if the pair of unitary operators $U$, $V$ is
  unitary equivalent to the pair of unitary operators $\tilde{U}$,
  $\tilde{V}$.

  Now, setting $S_{\tilde{U},\tilde{V}}=S_{U,V}$, rewrite the
  identity~(\ref{x1}) as follows:
  \begin{gather*}
    \End(S_{U,V})=\{\mathcal R\in B(\mathcal H)\,|\,\mathcal R(H_i)\subset
    H_i, i=1,\ldots,5\}
    \\
    \phantom{\End(S_{U,V})}{}= \{R\oplus R\in B(\mathcal H)\,|\,R\in
    B(H),RU=UR,RV=VR\}.
  \end{gather*}
  The latter identity immediately implies that the system $S_{U,V}$ is
  transitive if and only if the unitary operators $U$, $V$ are
  irreducible.
\end{proof}

Theorem~\ref{t1} allows to identify the description problem for
nonisomorphic transitive quintuples that correspond to f\/ive
orthogonal projections of a special type with that for
inequivalent irreducible pairs of unitary operators. The latter
problem is $*$-wild in the theory of $*$-representations of
$*$-algebras~\cite{KS,OS}.

\subsection[On $*$-wildness of the description problem for transitive
  systems corresponding to orthogonal projections with an additional
  relation]{On $\boldsymbol{*}$-wildness of the description problem for transitive
  systems\\ corresponding to orthogonal projections with an additional
  relation}

Let $P_1$, $P_2$, $P_3$ be orthogonal projections on a Hilbert
space $H$, and $P_2$, $P_3$ be mutually orthogo\-nal. Introduce a
system of f\/ive subspaces of the space $H$ corresponding to the
collection of orthogonal projections $P_1$, $P^\bot_1$,$P_2$,
$P_3$,$(P_2+P_3)^\bot$. Denote
\[
S_{P_1,P_2\bot P_3}=(H;\Imp P_1,\Imp P^\bot_1,\Imp P_2,\Imp
P_3,\Imp (P_2+P_3)^\bot).
\]

\begin{theorem}\label{t2}
  Let $P_1$, $P_2$, $P_3$ be orthogonal projections on a Hilbert space
  $H$ such that $P_2$ and~$P_3$ are mutually orthogonal, and
  $\tilde{P}_1$, $\tilde{P}_2$, $\tilde{P}_3$ be orthogonal
  projections on a Hilbert space $\tilde{H}$ such that~$\tilde{P}_2$
  and $\tilde{P}_3$ are mutually orthogonal. Then the system
  $S_{P_1,P_2\bot P_3}$ is transitive if and only if the projections
  $P_1$, $P_2$, $P_3$ are irreducible. Also, $S_{P_1,P_2\bot P_3}\cong
  S_{\tilde{P}_1,\tilde{P}_2\bot \tilde{P}_3}$ if and only if the
  triple of the orthogonal projections $P_1$, $P_2$, $P_3$ is unitary
  equivalent to the triple of the orthogonal projections
  $\tilde{P}_1$, $\tilde{P}_2$, $\tilde{P}_3$.
\end{theorem}

\begin{proof}
  Denote $H_1=\Imp P_1$, $H_2=\Imp P^\bot_1$, $H_3=\Imp P_2$,
  $H_4=\Imp P_3$, $H_5=\Imp (P_2+P_3)^\bot$, and let $\tilde{H}_1=\Imp
  \tilde{P}_1$, $\tilde{H}_2=\Imp \tilde{P}^\bot_1$, $\tilde{H}_3=\Imp
  \tilde{P}_2$, $\tilde{H}_4=\Imp \tilde{P}_3$, $\tilde{H}_5=\Imp
  (\tilde{P}_2+\tilde{P}_3)^\bot$.

  The proof of the theorem directly follows from the identity
  \begin{gather*}
  \{R\in B(H,\tilde{H})\,|\,
  R(H_i)\subset\tilde{H}_i,i=1,\ldots,5\}= \{R\in
  B(H,\tilde{H})\,|\,RP_i=\tilde{P}_iR,i=1,2,3\}.\!\!\!\! \tag*{\qed}
  \end{gather*} \renewcommand{\qed}{}
\end{proof}

Theorem~\ref{t2} identif\/ies the description problem for
nonisomorphic transitive quintuples of subspaces corresponding to
quintuples of orthogonal projections of a special type, the ones
such that their sum equals $2I_H$, with that for inequivalent
irreducible triples $P_1$, $P_2$, $P_3$ of ortho\-go\-nal
projections satisfying the condition $P_2\bot P_3$. The latter
problem is $*$-wild in the theory of $*$-representations of
$*$-algebras~\cite{KS,OS}.

\section[Transitive systems of subspaces corresponding to $\Rep\mathcal   P_{n,{\rm com}}$]{Transitive
systems of subspaces corresponding to $\boldsymbol{\Rep\mathcal
P_{n,{\rm com}}}$}

\subsection[On $*$-representations of the $*$-algebra $\mathcal P_{n,{\rm com}}$]{On $\boldsymbol{*}$-representations
of the $\boldsymbol{*}$-algebra $\boldsymbol{\mathcal P_{n,{\rm
com}}}$}

Denote by $\Sigma_n$ ($n\in\mathbb N$) the set $\alpha\in\mathbb
R_+$ such that there exists at least one $*$-representation of the
$*$-algebra $\mathcal P_{n,\alpha}=\mathbb C\langle
p_1,p_2,\ldots, p_n\,
|\,p^2_k=p^*_k=p_k,\sum\limits^n_{k=1}p_k=\alpha e\rangle$, i.e.,
the set of real numbers $\alpha$ such that there exist $n$
orthogonal projections $P_1, P_2, \ldots, P_n$ on a Hilbert space
$H$ satisfying $\sum\limits^n_{k=1}P_k=\alpha I_H$. It follows
from the def\/inition of the algebra $\mathcal P_{n,{\rm
com}}=\mathbb C\langle p_1,p_2,\ldots, p_n\,|\,p^2_k=p^*_k=p_k,
[\sum\limits_{k=1}^np_k,p_i]=0,\,\forall\, i=1,\ldots,n\rangle$
that all irreducible $*$-representations of $\mathcal P_{n,{\rm
com}}$ coincide with the union of irreducible $*$-representations
of $\mathcal P_{n,\alpha}$ taken over all $\alpha\in\Sigma_n$.

A description of the set $\Sigma_n$ for all $n\in\mathbb N$ is
obtained by S.A.~Kruglyak, V.I.~Rabanovich, and
Yu.S.~Samo\v\i{}lenko in~\cite{KRS}, and is given by
\begin{gather*}
\Sigma_2=\{0,1,2\},\quad\Sigma_3=\big\{0,1, \tfrac32,2,3\big\},
\\
\Sigma_n=\left\{\Lambda_n^0,\Lambda_n^1,
\left[\tfrac{n-\sqrt{n^2-4n}}2,\tfrac{n+\sqrt{n^2-4n}}2\right],
n-\Lambda_n^1,n-\Lambda_n^0\right\} \ \mbox{for}\ n\geq4,
\\
\Lambda_n^0=\Big\{
 0,1+ \tfrac1{n-1}, 1+
\tfrac1{(n-2)-\frac1{n-1}},\ldots,
1+\tfrac1{(n-2)-\tfrac1{(n-2)-\tfrac1{\ddots-\tfrac1{n-1}}}},\ldots
 \Big\},
\\
\Lambda_n^1=\Big\{
 1,1+ \tfrac1{n-2}, 1+
\tfrac1{(n-2)-\frac1{n-2}},\ldots, 1+
\tfrac1{(n-2)-\tfrac1{(n-2)-\tfrac1{\ddots-\tfrac1{n-2}}}},\ldots
\Big\}.
\end{gather*}
Here, the elements of the sets $\Lambda_n^0$, $\Lambda_n^1$,
$n-\Lambda_n^1$, $n-\Lambda_n^0$, in what follows, will be called
points of the discrete spectrum of the description problem for
unitary representations of the algebra~$\mathcal P_{n,{\rm com}}$,
whereas the elements of the line segment
$\left[\frac{n-\sqrt{n^2-4n}}2,\frac{n+\sqrt{n^2-4n}}2\right]$ are
called point of the conti\-nuous spectrum. For each point $\alpha$
in the sets $\Lambda_n^0$, $n-\Lambda_n^0$ there exists, up to
unitary equivalence, a unique irreducible $*$-representation of
the $*$-algebra $\mathcal P_{n,\alpha}$ and, hence, that
of~$\mathcal P_{n,{\rm com}}$. For each point $\alpha$ in the sets
$\Lambda_n^1$, $n-\Lambda_n^1$ there exist $n$ inequivalent
irreducible $*$-representations of the $*$-algebra $\mathcal
P_{n,\alpha}$ and, hence, those of $\mathcal P_{n,{\rm com}}$.

An important instrument for describing the set $\Sigma_n$ and
representations of $\mathcal P_{n,{\rm com}}$ is use of Coxeter
functors, constructed in~\cite{KRS}, between the categories of
$*$-representations of $\mathcal P_{n,\alpha}$ for dif\/ferent
values of the parameters.

Def\/ine a functor $\EuScript T:\Rep\,\mathcal
P_{n,\alpha}\rightarrow \Rep\,\mathcal P_{n,n-\alpha}$,
see~\cite{KRS}. Let $\pi$ be a representation of the
algebra~$\mathcal P_{n,\alpha}$, and $\pi(p_i)=P_i$,
$i=1,\ldots,n$, be orthogonal projections on a representation
space~$H$. Then the representation $\hat{\pi}=\EuScript T(\pi)$ in
$\Rep\mathcal P_{n,n-\alpha}$ is def\/ined by the identities
$\hat{\pi}(p_i)=(I-P_i)$ that give orthogonal projections on $H$.
We leave out a description of the action of the functor~$\EuScript
T$ on morphisms of the category $\Rep\mathcal P_{n,\alpha}$, since
it is not used in the sequel. Let us now def\/ine a~functor
$\EuScript S:\Rep\,\mathcal P_{n,\alpha}\rightarrow \Rep\,\mathcal
P_{n,\frac{\alpha}{\alpha-1}}$, see~\cite{KRS}. Again, let $\pi$
denote a representation in $\Rep\mathcal P_{n,\alpha}$, and by
$P_1, P_2, \ldots, P_n$ denote the corresponding orthogonal
projections on the representation space $H$. Consider the
subspaces $H_i=\Imp P_i$ ($i=1,\ldots,n$).  Let
$\Gamma_i:H_i\rightarrow H$, $i=1,\ldots,n$, be the natural
isometries. Then
\begin{gather}\label{x2}
  \Gamma_i^*\Gamma_i=I_{H_i},\qquad \Gamma_i\Gamma_i^*=P_i,\qquad
  i=1,\ldots,n.
\end{gather}

Let an operator $\Gamma$ be def\/ined by the matrix
$\Gamma=[\Gamma_1,\Gamma_2,\ldots,\Gamma_n]:\EuScript H=H_1\oplus
H_2\oplus\cdots \oplus H_n\rightarrow H$. Then the natural
isometry $\sqrt{\frac{\alpha-1}{\alpha}}\Delta^*$ that acts from
the orthogonal complement $\hat{H}$ to the subspace $\Imp
\Gamma^*$ into the space~$\mathcal H$ def\/ines the isometries
$\Delta_k=\Delta|_{\Imp
  P_k}:H_k\rightarrow\hat{H}$, $k=1,\ldots,n$. The orthogonal
projections $Q_i=\Delta_i\Delta_i^*$, $i=1,\ldots,n$, on the space
$\hat{H}$ make the corresponding representation in $\EuScript
S(\Rep\,\mathcal P_{n,\alpha})$, i.e.\ the representation
$\hat{\pi}=\EuScript S(\pi)$ in $\Rep\mathcal
P_{n,\frac{\alpha}{\alpha-1}}$ is given by the identities
$\hat{\pi}(p_i)=Q_i$. Write down the relations satisf\/ied by the
operators $\{\Delta_i\}^n_i$,
\begin{gather}\label{x3}
\Delta^*_i\Delta_i=I_{H_i},\qquad \Delta_i\Delta_i^*=Q_i,\qquad
i=1,\ldots,n.
\end{gather}

We will not describe the action of the functor $\EuScript S$ on
morphisms of the category $\Rep\mathcal P_{n,\alpha}$, since we
will not use it in the sequel.

Following~\cite{KRS}, introduce a functor $\Phi^+:\Rep\mathcal
P_{n,\alpha}\rightarrow\Rep\mathcal P_{n,1+\frac1{n-1-\alpha}}$
def\/ined by $\Phi^+(\pi)=\EuScript S(\EuScript T(\pi))$ for
$\alpha<n-1$. Denote by $\pi_k$ ($k=0,1,\ldots,n$) the following
representations in $\Rep\mathcal P_{n,\alpha}$: $\pi_0(p_i)=0$,
$i=1,\ldots,n$, where the space of representation is $\mathbb C$;
$\pi_k(p_i)=0$ if $i\neq k$ and $\pi_k(p_k)=1$, $k=1,\ldots,n$,
with $\mathbb C$ as the representation space. For an arbitrary
irreducible representation $\pi$ of the algebra $\mathcal
P_{n,\alpha}$ in the case of points of the discrete spectrum, one
can assert that either $\pi$ or $\EuScript T(\pi)$ is unitary
equivalent to a representation of the form
$\Phi^{+s}(\check{\pi})$, where the representation $\check{\pi}$
is one of the simplest representations $\pi_k$,
$k=\overline{0,n}$, and $s$ is a natural number.

\subsection[Transitive systems of $n$ subspaces corresponding to
  $\Rep\mathcal P_{n,\alpha}$]{Transitive systems of $\boldsymbol{n}$ subspaces corresponding to
  $\boldsymbol{\Rep\mathcal P_{n,\alpha}}$}

The systems of subspaces, $S_{\pi_k}$, $k=0,1,\ldots,n$, are
clearly nonisomorphic transitive systems of~$n$ subspaces of the
space $\mathbb C$. I.M.~Gel'fand and V.A.~Ponamarev in~\cite{GP},
by using the functor technique, construct from the systems
$S_{\pi_k}$, $k=0,1,\ldots,n$, inf\/inite series of indecomposable
systems, which turn out to be are transitive, of $n$ subspaces. In
this section we show that the Coxeter functors in~\cite{KRS}, as
the functors in~\cite{GP}, transform nonisomorphic transitive
systems into nonisomorphic transitive systems and, consequently,
all systems of the form $S_{\Phi^{+s}(\check{\pi})}$ and
$S^\bot_{\Phi^{+s}(\check{\pi})}$, where the representation
$\check{\pi}$ is one of the simplest representations $\pi_k$,
$k=0,1,\ldots,n$, and $s$ is a natural number, will be
nonisomorphic transitive systems. Hence, we have the following
theorem.

\begin{theorem}\label{t3}
  Systems of $n$ subspaces $S_{\pi}$ constructed from irreducible
  inequivalent representations $\pi\in\Rep\mathcal P_{n,\alpha}$, for
  $\alpha$ in the discrete spectrum, are nonisomorphic and
  transitive.
\end{theorem}

To prove the theorem, by using the Coxeter functors $\EuScript T$
and $\EuScript S$ in~\cite{KRS}, we construct auxiliary functors
$\EuScript T'$ and $\EuScript S'$. The action of the functors
$\EuScript T':\Rep\,\mathcal P_{n,\alpha}\rightarrow
\Rep\,\mathcal P_{n,n-\alpha}$ and $\EuScript S':\Rep\,\mathcal
P_{n,\alpha}\rightarrow \Rep\,\mathcal
P_{n,\frac{\alpha}{\alpha-1}}$ on the objects of the category is
def\/ined to coincide with the actions of $\EuScript T$ and
$\EuScript S$, that is, $\EuScript T(\pi)=\EuScript T'(\pi)$ and
$\EuScript S(\pi)=\EuScript S'(\pi)$ $\forall\, \pi\in\Rep\mathcal
P_{n,\alpha}$. The morphisms of the category of representations
are def\/ined dif\/ferently. Let $\pi\in\Rep(\mathcal
P_{n,\alpha},H)$ and $\tilde{\pi}\in\Rep(\mathcal
P_{n,\alpha},\tilde{H})$. A linear operator $C\in B(H,\tilde{H})$
is called a morphism of the category of representations,
$C\in\Mor(\pi,\tilde{\pi})$, if
$C\pi(p_i)=\tilde{\pi}(p_i)C\pi(p_i)$, $i=1,\ldots,n$, that is,
\begin{gather}\label{x4}
  CP_i=\tilde{P}_iCP_i, \qquad i=1,\ldots,n.
\end{gather}

The restrictions $C|_{H_i}$, $i=1,\ldots,n$, are denoted by $C_i$.
Let us show that the operators $C_i$ map $H_i$ into $\tilde{H}_i$,
that is,
\begin{gather}\label{x5}
  C_i(H_i)\subset\tilde{H}_i,\qquad i=1,\ldots,n.
\end{gather}
Indeed, for $x\in H_i$, we have $C_ix=Cx=CP_ix=\tilde{P}_iCP_ix$
and, consequently, $C_ix\in\tilde{H}_i$.

If $x\in H_i$, then~(\ref{x4}) and~(\ref{x5}) give
\[
C\Gamma_ix=Cx=CP_ix=\tilde{P}_iCP_ix=
\tilde{P}_iC_ix=C_ix=\tilde{\Gamma}_iC_ix,
\]
so that
\begin{gather}\label{x6}
  C\Gamma_i=\tilde{\Gamma}_iC_i,\qquad i=1,\ldots,n.
\end{gather}

The identities~(\ref{x4}) are equivalent to the inclusions
$C(H_i)\subset\tilde{H}_i$, $i=1,\ldots,n$, which imme\-diately
gives the following relations:
\begin{gather}\label{x7}
  C_i=\tilde{\Gamma}_i^*C\Gamma_i,\qquad i=1,\ldots,n.
\end{gather}

Formula~(\ref{x6}) allows to represent $C$ as
\begin{gather}\label{x8}
  C=\frac1{\alpha}\sum_{i=1}^n\tilde{\Gamma}_iC_i\Gamma_i^*.
\end{gather}

Indeed,
$\frac1{\alpha}\sum\limits_{i=1}^n\tilde{\Gamma}_iC_i\Gamma_i^*=
\frac1{\alpha}\sum\limits_{i=1}^nC
\Gamma_i\Gamma_i^*=C(\frac1{\alpha} \sum\limits_{i=1}^nP_i)=C$.

Consider an operator $\hat{C}:\hat{H}\rightarrow\hat{\tilde{H}}$
def\/ined by
\begin{gather}\label{x9}
  \hat{C}=\frac{\alpha-1}{\alpha}\sum_{i=1}^n\tilde{\Delta}_iC_i\Delta^*_i.
\end{gather}
Using the following properties of the operators~\cite{KRS}
$\{\Gamma_i\}_{i=1}^n$, $\{\Gamma_i^*\}_{i=1}^n$,
$\{\Delta_i\}_{i=1}^n$, $\{\Delta_i^*\}_{i=1}^n$:
\begin{gather}\label{x10}
\sum_{i=1}^n\Gamma_i\Delta_i^*=0,
\\
\label{x11}
\Delta_i^*\Delta_j=-\frac1{\alpha-1}\Gamma_i^*\Gamma_j,\quad\;
i\neq j,
\end{gather}
let us prove that
\begin{gather}\label{x12}
  \tilde{\Delta}^*_k\hat{C}=C_k\Delta^*_k,\quad  k=1,\ldots,n.
\end{gather}
Indeed,
\begin{gather*}
\tilde{\Delta}^*_k\hat{C}=
\tilde{\Delta}^*_k\left(\frac{\alpha-1}{\alpha}\sum_{i=1}^n\tilde{\Delta}_iC_i
\Delta^*_i\right)=
\frac{\alpha-1}{\alpha}\sum_{i=1}^n(\tilde{\Delta}^*_k\tilde{\Delta}_i)C_i
\Delta^*_i\\
\phantom{\tilde{\Delta}^*_k\hat{C}}{} =
\frac{\alpha-1}{\alpha}(\tilde{\Delta}^*_k\tilde{\Delta}_k)C_k\Delta^*_k+
\frac{\alpha-1}{\alpha}\sum_{i=1\atop i\neq
  k}^n(-\frac1{\alpha-1})\tilde{\Gamma}^*_k(\tilde{\Gamma}_iC_i)\Delta^*_i\\
\phantom{\tilde{\Delta}^*_k\hat{C}}{}  =
\frac{\alpha-1}{\alpha}C_k\Delta^*_k- \frac1{\alpha}\sum_{i=1\atop
  i\neq k}^n\tilde{\Gamma}^*_k(\tilde{\Gamma}_iC_i)\Delta^*_i=
\frac{\alpha-1}{\alpha}C_k\Delta^*_k-
\frac1{\alpha}\tilde{\Gamma}^*_kC\sum_{i=1\atop i\neq
  k}^n\Gamma_i\Delta^*_i\\
\phantom{\tilde{\Delta}^*_k\hat{C}}{}  =
\frac{\alpha-1}{\alpha}C_k\Delta^*_k-
\frac1{\alpha}\tilde{\Gamma}^*_kC\left(\sum_{i=1}^n\Gamma_i\Delta^*_i-\Gamma_k
\Delta^*_k\right)= \frac{\alpha-1}{\alpha}C_k\Delta^*_k+
\frac1{\alpha}(\tilde{\Gamma}^*_k)C\Gamma_k\Delta^*_k=C_k\Delta^*_k.\!\!
\end{gather*}

Now, let us show that
\begin{gather}\label{x13}
  C_k=\tilde{\Delta}^*_k\hat{C}\Delta_k,\qquad k=1,\ldots,n.
\end{gather}

Using~(\ref{x2}), (\ref{x3}), (\ref{x7}), (\ref{x8}), (\ref{x9}),
and (\ref{x11}) we get
\begin{gather*}
\tilde{\Delta}_k^*\hat{C}\Delta_k=\tilde{\Delta}_k^*
\left(\frac{\alpha-1}{\alpha}\sum_{i=1}^n\tilde{\Delta}_iC_i\Delta^*_i\right)
\Delta_k=
\frac{\alpha-1}{\alpha}\sum_{i=1}^n\tilde{\Delta}_k^*\tilde{\Delta}_iC_i
\Delta^*_i\Delta_k\\
\phantom{\tilde{\Delta}_k^*\hat{C}\Delta_k}{}
=\frac{\alpha-1}{\alpha}\tilde{\Delta}_k^*\tilde{\Delta}_kC_k
\Delta^*_k\Delta_k+ \frac{\alpha-1}{\alpha}\sum_{i=1\atop i\neq
  k}^n\tilde{\Delta}_k^*\tilde{\Delta}_iC_i\Delta^*_i\Delta_k\\
\phantom{\tilde{\Delta}_k^*\hat{C}\Delta_k}{}=
\frac{\alpha-1}{\alpha}C_k+\frac1{\alpha(\alpha-1)}\sum_{i=1\atop
  i\neq k}^n\tilde{\Gamma}_k^*\tilde{\Gamma}_iC_i\Gamma^*_i\Gamma_k\\
\phantom{\tilde{\Delta}_k^*\hat{C}\Delta_k}{}  =
\frac{\alpha-1}{\alpha}C_k+\frac1{\alpha-1}\tilde{\Gamma}_k^*
\left(\frac1{\alpha}\sum_{i=1\atop i\neq
  k}^n\tilde{\Gamma}_iC_i\Gamma_i^*\right)\Gamma_k\\
\phantom{\tilde{\Delta}_k^*\hat{C}\Delta_k}{}  =
\frac{\alpha-1}{\alpha}C_k+\frac1{\alpha-1}\tilde{\Gamma}_k^*C\Gamma_k
-\frac1{\alpha(\alpha-1)}\tilde{\Gamma}_k^*\tilde{\Gamma}_kC_k
\Gamma^*_k\Gamma_k =C_k.
\end{gather*}

Now, it follows from~(\ref{x12}) and~(\ref{x13}) that
$\tilde{Q}_k\hat{C}=\tilde{\Delta}_k\tilde{\Delta}^*_k\hat{C}=
\tilde{\Delta}_kC_k\Delta^*_k=\tilde{Q}_k\hat{C}Q_k$, that is,
$\tilde{Q}_k\hat{C}=\tilde{Q}_k\hat{C}Q_k$, $k=1,\ldots,n$.
Whence,
\begin{gather}\label{x14}
  \hat{C}^*\tilde{Q}_k=Q_k\hat{C}^*\tilde{Q}_k, \qquad
  k=1,\ldots,n.
\end{gather}

The latter identities mean that $\hat{C}^*\in\Mor(\EuScript
S'(\tilde{\pi}),\EuScript S'(\pi))$. The action of the auxiliary
functors~$\EuScript T'$ and $\EuScript S'$ on morphisms of the
category $\Rep\mathcal P_{n,\alpha}$ are def\/ined by $\EuScript
T'(C)=C^*$ and $\EuScript S'(C)=\hat{C}^*$ for any
$C\in\Mor(\pi,\tilde{\pi})$. This completes the construction of
the auxiliary functors.

\begin{lemma}\label{lem1}
  The functors $\EuScript T'$ and $\EuScript S'$ are category
  equivalences.
\end{lemma}

\begin{proof}
  It is easy to check by using the def\/inition that the functor
  $\EuScript T'$ is univalent and complete. $\EuScript T^2=\Id$ and
  $\EuScript T'(\pi)=\EuScript T(\pi)$ for any $\pi\in\Rep\mathcal
  P_{n,\alpha}$.  Consequently, the functor $\EuScript T'$ is an
  equivalence between the categories $\Rep\mathcal P_{n,\alpha}$ and
  $\Rep\mathcal P_{n,n-\alpha}$.

  Now, let us prove the lemma for the functor $\EuScript S'$. Let us
  show that the functor $\EuScript S'$ is univalent. Let
  $C,D\in\Mor(\pi,\tilde{\pi})$ and $C\neq D$, and show that
  $\EuScript S'(C)\neq \EuScript S'(D)$. Indeed, if $\EuScript
  S'(C)= \EuScript S'(D)$, then $\hat{C}^*=\hat{D}^*$ and
  $\hat{C}=\hat{D}$. By~(\ref{x13}), we have
  \[
  C_i=\tilde{\Delta}^*_i\hat{C}\Delta_i=\tilde{\Delta}^*_i\hat{D}
  \Delta_i=D_i,\qquad i=1,\ldots,n.
  \]
 {\samepage Using the decomposition~(\ref{x8}) we get
  \[
  C=\frac1{\alpha}\sum_{i=1}^n\tilde{\Gamma}_iC_i\Gamma_i^*,\qquad
  D=\frac1{\alpha}\sum_{i=1}^n\tilde{\Gamma}_iD_i\Gamma_i^*.
  \]
  Then $C=D$ and, hence, the functor $\EuScript S'$ is univalent.}

  Let us now show that $\EuScript S'$ is complete. Let
  $R\in\Mor(\EuScript S'(\tilde{\pi}),\EuScript S'(\pi))$. To prove
  the completeness, construct a linear operator from the set
  $\Mor(\pi,\tilde{\pi})$ such that the functor takes on this morphism
  the value $R$. Since $R\in\Mor(\EuScript S'(\tilde{\pi}),\EuScript
  S'(\pi))$, the operator $R:\hat{\tilde{H}}\rightarrow\hat{H}$
  satisf\/ies
  \[
  R\tilde{Q}_k=Q_kR\tilde{Q}_k,\qquad k=1,\ldots,n.
  \]

  Consider an operator $\hat{r}$ in $B(\hat{H},\hat{\tilde{H}})$
  such that $\hat{r}^*=R$. Then the former identities can be written
  as
  \[
  \hat{r}^*\tilde{Q}_k=Q_k\hat{r}^*\tilde{Q}_k,\qquad
  k=1,\ldots,n,
  \]
  and, consequently,
  \begin{gather}\label{x15}
    \tilde{Q}_k\hat{r}=\tilde{Q}_k\hat{r}Q_k,\qquad k=1,\ldots,n.
  \end{gather}

  Denote by $r_k$ the operators
  $r_k=\tilde{\Delta}^*_k\hat{r}\Delta_k:H_k\rightarrow\tilde{H}_k$,
  $k=1,\ldots,n$, and show that $\hat{r}$ can be represented as
  \begin{gather}\label{x16}
    \hat{r}=\frac{\alpha-1}{\alpha}\sum_{k=1}^n
    \tilde{\Delta}_kr_k\Delta_k^*.
  \end{gather}
  Indeed,
  \[
  \frac{\alpha-1}{\alpha}\sum_{k=1}^n\tilde{\Delta}_kr_k\Delta_k^*=
  \frac{\alpha-1}{\alpha}\sum_{k=1}^n\tilde{\Delta}_k\tilde{\Delta}^*_k
  \hat{r}\Delta_k
  \Delta_k^*=\frac{\alpha-1}{\alpha}\sum_{k=1}^n\tilde{Q}_k\hat{r}Q_k=
  \left(\frac{\alpha-1}{\alpha}\sum_{k=1}^n\tilde{Q}_k\right)\hat{r}=\hat{r}.
  \]

  It follows from the def\/inition of $r_k$ and identities~(\ref{x3}),
  (\ref{x15}) that
  \begin{gather*}
    r_k\Delta_k^*=(\tilde{\Delta}^*_k\hat{r}\Delta_k)\Delta_k^*=
  \tilde{\Delta}^*_k\hat{r}(\Delta_k\Delta_k^*)=
  \tilde{\Delta}^*_k\hat{r}Q_k=
  I_{\tilde{H}_k}\tilde{\Delta}^*_k\hat{r}Q_k=
  (\tilde{\Delta}^*_k\tilde{\Delta}_k)
  \tilde{\Delta}^*_k\hat{r}Q_k\\
 \phantom{r_k\Delta_k^*}{} =\tilde{\Delta}^*_k(\tilde{\Delta}_k
  \tilde{\Delta}^*_k)\hat{r}Q_k=\tilde{\Delta}^*_k\tilde{Q}_k\hat{r}Q_k=
  \tilde{\Delta}^*_k\tilde{Q}_k\hat{r}=\tilde{\Delta}^*_k(\tilde{\Delta}_k
  \tilde{\Delta}^*_k)\hat{r}=(\tilde{\Delta}^*_k\tilde{\Delta}_k)
  \tilde{\Delta}^*_k\hat{r}=\tilde{\Delta}^*_k\hat{r}.
\end{gather*}
Hence, we have
  \begin{gather}\label{x17}
    r_k\Delta_k^*=\tilde{\Delta}^*_k\hat{r},\qquad k=1,\ldots,n.
  \end{gather}
  Consider the operator
  \begin{gather}\label{x18}
    r=\frac1{\alpha}\sum_{i=1}^n\tilde{\Gamma}_ir_i\Gamma^*_i.
  \end{gather}
  Using~(\ref{x2}), (\ref{x10}), (\ref{x11}), (\ref{x17}) we get
  \begin{gather}\label{x19}
    r\Gamma_k=\tilde{\Gamma}_kr_k,\qquad k=1,\ldots,n,
  \\
  \label{x20}
    r_k=\tilde{\Gamma}_k^*r\Gamma_k,\qquad k=1,\ldots,n.
  \end{gather}
  Indeed,
  \begin{gather*}
    r\Gamma_k=\frac1{\alpha}\sum_{i=1}^n\tilde{\Gamma}_ir_i\Gamma_i^*\Gamma_k=
  \frac1{\alpha}\tilde{\Gamma_k}r_k+\frac1{\alpha}\sum_{i=1\atop i\neq
    j }^n\tilde{\Gamma}_ir_i(\Gamma_i^*\Gamma_k)
  =\frac1{\alpha}\tilde{\Gamma}_kr_k-\frac{\alpha-1}{\alpha}\sum_{i=1\atop
    i\neq j }^n\tilde{\Gamma}_i(r_i\Delta_i^*)\Delta_k\\
  \phantom{r\Gamma_k}{}=
  \frac1{\alpha}\tilde{\Gamma}_kr_k-\frac{\alpha-1}{\alpha}\sum_{i=1\atop
    i\neq j }^n\tilde{\Gamma}_i(\tilde{\Delta}_i^*\hat{r})
  \Delta_k=\frac1{\alpha}\tilde{\Gamma}_kr_k
  +\frac{\alpha-1}{\alpha}\tilde{\Gamma}_k\Delta_k^*\hat{r}
  \Delta_k=\tilde{\Gamma}_kr_k
\end{gather*}
and
\begin{gather*}
  \tilde{\Gamma}_k^*r\Gamma_k=
  \frac1{\alpha}\tilde{\Gamma}_k^*\left(\sum_{i=1}^n\tilde{\Gamma}_i
  r_i\Gamma_i^*\right)\Gamma_k=
  \frac1{\alpha}r_k+\frac1{\alpha}\sum_{i=1\atop i\neq
    j}^n\tilde{\Gamma}_k^*\tilde{\Gamma}_ir_i\Gamma_i^*\Gamma_k\\
  \phantom{\tilde{\Gamma}_k^*r\Gamma_k}{} =
  \frac1{\alpha}r_k+\frac{(\alpha-1)^2}{\alpha}\sum_{i=1\atop i\neq
    j}^n\tilde{\Delta}_k^*\tilde{\Delta}_ir_i
  \Delta_i^*\Delta_k
  =\frac1{\alpha}r_k+
  (\alpha-1)\tilde{\Delta}_k^*\hat{r}\Delta_k
  -\frac{(\alpha-1)^2}{\alpha}r_k=r_k.
\end{gather*}
It follows from~(\ref{x19})
  and (\ref{x20}) that
  $rP_k=r\Gamma_k\Gamma_k^*=\tilde{\Gamma}_kr_k\Gamma_k^*=
  \tilde{\Gamma}_k\tilde{\Gamma}_k^*r_k\Gamma_k\Gamma_k^*=
  \tilde{P}_krP_k$, which means that $r\in\Mor(\pi,\tilde{\pi})$.

  It is easy to check that $\EuScript S'(r)=R$ and, consequently, the
  functor $\EuScript S'$ is complete. So the univalence and
  completeness properties of the functor $\EuScript S'$ are checked,
  $\EuScript S^2=\Id$ and $\EuScript S'(\pi)=\EuScript S(\pi)$ for any
  $\pi\in\Rep\mathcal P_{n,\alpha}$. Consequently, the functor
  $\EuScript S'$ is an equivalence between the categories
  $\Rep\mathcal P_{n,\alpha}$ and $\Rep\mathcal
  P_{n,\frac{\alpha}{\alpha-1}}$.
\end{proof}

\begin{lemma}\label{lem2}
  If a system $S_{\pi}$, $\pi\in\Rep\mathcal P_{n,{\rm com}}$, of subspaces
  is transitive, then the system $S_{\Phi^+(\pi)}$ of subspaces is
  transitive. Here, $S_{\pi}\cong S_{\tilde{\pi}}$ if and only if
  $S_{\Phi^+(\pi)}\cong S_{\Phi^+(\tilde{\pi})}$.
\end{lemma}

\begin{proof}
  For the functors $\EuScript T$ and $\EuScript S$, we have $\EuScript
  T(\pi)=\EuScript T'(\pi)$ and $\EuScript S(\pi)=\EuScript S'(\pi)$
  for any $\pi\in\Rep\mathcal P_{n,\alpha}$. Consequently,
  $S_{\EuScript T(\pi)}=S_{\EuScript T'(\pi)}$ and $S_{\EuScript
    S(\pi)}=S_{\EuScript S'(\pi)}$. By Lemma~\ref{lem1}, $\EuScript
  T'$ is an equivalence of the categories that shows that if a
  system $S_{\pi}$, $\pi\in\Rep\mathcal P_{n,{\rm com}}$, of subspaces is
  transitive, then the system $S_{\EuScript T(\pi)}$ of subspaces is
  transitive. We also have that $S_{\pi}\cong S_{\tilde{\pi}}$ if and
  only if $S_{\EuScript T(\pi)}\cong S_{\EuScript T(\tilde{\pi})}$.

  Let us now consider the systems $S_{\EuScript S(\pi)}$,
  $\pi\in\Rep\mathcal P_{n,\alpha}$ of subspaces. Let
  $\pi,\tilde{\pi}\in \Rep\mathcal P_{n,\alpha}$. Consider the systems
  of subspaces $S_{\pi}=(H;H_1,H_2,\ldots,H_n)$ and
  $S_{\tilde{\pi}}=(\tilde{H};\tilde{H}_1,\tilde{H}_2,\ldots,\tilde{H}_n)$,
  that, respectively, correspond to the representations $\pi$ and
  $\tilde{\pi}$. Let the systems of subspaces be isomorphic, that is,
  $S_{\pi}\cong S_{\tilde{\pi}}$. By the def\/inition of isomorphic
  systems, there exists a linear operator $T\in B(H,\tilde{H})$ such
  that $T^{-1}\in B(\tilde{H},H)$ and $T(H_i)=\tilde{H}_i$,
  $i=1,\ldots,n$. It follows from $T(H_i)=\tilde{H}_i$,
  $i=1,\ldots,n$, that $T(H_i)\subset\tilde{H}_i$,
  $i=1,\ldots,n$, and, consequently, we get the relations
  $TP_i=\tilde{P}_iTP_i$, $i=1,\ldots,n$. The latter relations
  mean that $T\in\Mor(\pi,\tilde{\pi})$ if $\hat{T}^*\in\Mor(\EuScript
  S'(\tilde{\pi}),\EuScript S'(\pi))$, and
  \begin{gather}\label{x21}
    \hat{T}^*(\Imp \tilde{Q}_i)\subset(\Imp Q_i),\qquad
    i=1,\ldots,n.
  \end{gather}

  Again, using $T(H_i)=\tilde{H}_i$, $i=1,\ldots,n$, we get
  $T(H_i)\supset\tilde{H}_i$, $i=1,\ldots,n$, so that
  $T^{-1}(\tilde{H}_i)\subset H_i$, $i=1,\ldots,n$, and,
  respectively, $T^{-1}\tilde{P}_i=P_iT^{-1}\tilde{P}_i$,
  $i=1,\ldots,n$. This means that
  $T^{-1}\in\Mor(\tilde{\pi},\pi)$, hence,
  $\widehat{T^{-1}}^*\in\Mor(\EuScript S'(\pi),\EuScript
  S'(\tilde{\pi}))$, and using
  $\widehat{T^{-1}}^*=(\hat{T}^{-1})^*=(\hat{T}^*)^{-1}$ we get
  \[
  \Imp \tilde{Q}_i\supset(\hat{T}^*)^{-1}(\Imp Q_i),\qquad
  i=1,\ldots,n,
  \]
  so that
  \begin{gather}\label{x22}
    \hat{T}^*(\Imp \tilde{Q}_i)\supset \Imp Q_i ,\qquad
    i=1,\ldots,n.
  \end{gather}
  It follows from~(\ref{x21}) and (\ref{x22}) that
  \[
  \hat{T}^*(\Imp
  \tilde{Q}_i)= \Imp Q_i ,\qquad i=1,\ldots,n,
  \]
  i.e., it is an isomorphism of the systems corresponding to the
  representations $\EuScript S'(\pi)$ and $\EuScript S'(\tilde{\pi})$
  and, since the functors $\EuScript S'$ and $\EuScript S$ coincide on
  the objects of the categories, it is an isomorphism of the systems
  corresponding to the representations $\EuScript S(\pi)$ and
  $\EuScript S(\tilde{\pi})$.

  Since $\EuScript S'$ is complete, using similar reasonings it is
  easy to show that the functor $\EuScript S'$ and, hence, $\EuScript
  S$ takes the representations corresponding to nonisomorphic systems
  to representations that also correspond to nonisomorphic systems.

  Let again $\pi$ be a representation of the algebra $\mathcal
  P_{n,\alpha}$, and $\pi(p_i)=P_i$, $i=1,\ldots,n$, be orthogonal
  projections on a representation space $H$. Assume that the system of
  projections $P_1, P_2, \ldots, P_n$ gives rise to a transitive
  system of subspaces $S_{\pi}=(H;H_1,H_2,\ldots,H_n)$, where
  $H_i=P_iH$, $i=1,\ldots,n$, that is,
  \[
  \End(S_{\pi})=\{r\in B(H)\,|\,r(H_i)\subset H_i,
  i=1,\ldots,n\}=\Mor(\pi,\pi)=\mathbb CI.
  \]

  Consider $\EuScript S'(\pi)=\hat{\pi}$, where $\hat{\pi}(q_i)=Q_i$,
  $i=1,\ldots,n$, and the corresponding system of subspaces
  $S_{\hat{\pi}}$. Let now $R\in \End(S_{\hat{\pi}})$. Since
  $\End(S_{\hat{\pi}})=\Mor(\EuScript S'(\pi),\EuScript S'(\pi))$ and
  the functor $\EuScript S'$ is complete, we see that $\EuScript
  S'(r)=R$, where $r\in\Mor(\pi,\pi)$ is constructed from the operator
  $R^*=\frac{\alpha-1}{\alpha}\sum_{k=1}^n\tilde{\Delta}_kr_k\Delta_k^*$,
  $r_i=\tilde{\Delta}^*_iR^*\Delta_i:H_i\rightarrow\tilde{H}_i$,
  $i=1,\ldots,n$, as follows:
  \begin{gather}\label{d15}
    r=\frac1{\alpha}\sum_{i=1}^n\Gamma_ir_i\Gamma^*_i.
  \end{gather}
  By using $R\in \Mor(\EuScript S'(\pi),\EuScript S'(\pi))$, we
  obtain, similarly to~(\ref{x20}), that
  \begin{gather}\label{x23}
    r_i=\Gamma_i^*r\Gamma_i,\qquad i=1,\ldots,n.
  \end{gather}

  Since the system $S_{\pi}$ is transitive, the operator $r$ is a
  scalar, that is, $r=\lambda I_{H}$. Using that
  $\Gamma_i^*\Gamma_i=I_{H_i}$, $i=1,\ldots,n$, and~(\ref{x23}) we
  get
  \[
  r_i=\lambda I_{H_i},\qquad i=1,\ldots,n.
  \]
  Then $R^*$ is a scalar operator and, consequently, $R$ is also a
  scalar operator that means that the system $S_{\EuScript S'(\pi)}$
  is transitive and such is $S_{\EuScript S(\pi)}$.
\end{proof}

The statement of Theorem~\ref{t3} follows directly from
Lemma~\ref{lem2}.

\section[Transitive systems of subspaces corresponding to $\Rep\mathcal  P_{n,{\rm abo},\tau}$]{Transitive
 systems of subspaces corresponding to $\boldsymbol{\Rep\mathcal  P_{n,{\rm abo},\tau}}$}

\subsection[Equivalence of the categories $\Rep\mathcal P_{n,\alpha}$
  and $\Rep\mathcal P_{n,{\rm abo},\tau}$]{Equivalence of the categories $\boldsymbol{\Rep\mathcal P_{n,\alpha}}$
  and $\boldsymbol{\Rep\mathcal P_{n,{\rm abo},\tau}}$}

Let us examine the equivalence $\EuScript F$, constructed in
\cite{PS}, between the categories of
$*$-representa\-tions~$\mathcal P_{n,\alpha}$ and $\mathcal
P_{n,{\rm abo},\frac1{\alpha}}$, $\alpha\neq 0$. Theorem~\ref{t3} allows
to consider nonisomorphic transitive systems of $n$ subspaces of
the form $S_{\pi}$, constructed from representations of the
algebras $\mathcal P_{n,\alpha}$ for $\alpha$ lying in the
discrete spectrum. The equivalence $\EuScript F$, in its turn,
allows to construct nonisomorphic transitive systems $S_{\EuScript
F(\pi)}$ of $n+1$ subspaces starting with nonisomorphic transitive
systems $S_{\pi}$, $\pi\in\mathcal P_{n,\alpha}$, of $n$
subspaces.

Let us describe the equivalence $\EuScript F$. Let $\pi$ be a
representation of the algebra $\mathcal P_{n,\alpha}$, and
$\pi(p_i)=P_i$, $i=1,\ldots,n$, be orthogonal projections on a
representation space $H$. As it was done in Section~4, let us
introduce the spaces $H_i=\Imp P_i$ and the natural isometries
$\Gamma_i:H_i\rightarrow H$. Let $\mathcal H=H_1\oplus
H_2\oplus\cdots\oplus H_n$. Def\/ine a linear operator $\Gamma
:H_1\oplus H_2\oplus\cdots\oplus H_n\rightarrow H$ in terms of the
matrix $\Gamma=(\Gamma_1\;\Gamma_2\;\ldots\Gamma_n)$ of the
dimension $n\times 1$. Let $Q_i$ denote $n$ orthogonal
projections, $Q_i=\Diag(0,\ldots,0,I_{H_i},0,\ldots,0)$,
$i=1,\ldots,n$, and $P:\mathcal H\rightarrow\mathcal H$ an
orthogonal projection def\/ined by
$P=\frac1{\alpha}\Gamma^*\Gamma$ with the block matrix
$P=\frac1{\alpha}||\Gamma_i^*\Gamma_j||_{i,j=1}^n$ on the space
$H_1\oplus H_2\oplus\cdots\oplus H_n$.

Let a functor $\EuScript F:\Rep\mathcal
P_{n,\alpha}\rightarrow\Rep\mathcal P_{n,{\rm abo},\frac1{\alpha}}$,
$\alpha\neq 0$, be def\/ined on objects of the category of
representations as follows: $\EuScript F(\pi)=\hat{\pi}$, where
$\hat{\pi}(q_i)=Q_i$, $i=1,\ldots,n$, and $\hat{\pi}(p)=P$. The
identities $\sum\limits_{i=1}^nQ_i=I$ and
$Q_iPQ_i=\frac1{\alpha}Q_i$, $i=1,\ldots,n$, are easily checked.
We do not describe the action of the functor $\EuScript F$ on
morphisms of the category $\Rep\mathcal P_{n,\alpha}$, since we
will not use it.

\begin{theorem}\label{t4}
  Systems of $n+1$ subspaces,  $S_{\EuScript F(\pi)}$, constructed from
  irreducible inequivalent representations $\pi\in\Rep\mathcal
  P_{n,\alpha}$, where $\alpha$ is in the discrete spectrum, are
  nonisomorphic and transitive.
\end{theorem}

To prove the theorem, construct an auxiliary functor $\EuScript
F':\Rep\mathcal P_{n,\alpha}\rightarrow\Rep\mathcal
P_{n,{\rm abo},\frac1{\alpha}}$, $\alpha\neq 0$, the action of which on
objects coincides with the action of $\EuScript F$, that is,
$\EuScript F'(\pi)=\EuScript F(\pi)$ for all $\pi\in\Rep\mathcal
P_{n,\alpha}$. Morphisms are def\/ined as in Section~4. Let
$\pi\in\Rep(\mathcal P_{n,\alpha},H)$ and
$\tilde{\pi}\in\Rep(\mathcal P_{n,\alpha},\tilde{H})$. A linear
operator $C\in B(H,\tilde{H})$ will be called a morphism of the
category of representations, written $C\in\Mor(\pi,\tilde{\pi})$,
if $C\pi(p_i)=\tilde{\pi}(p_i)C\pi(p_i)$, that is,
  \begin{gather}\label{x25}
    CP_i=\tilde{P}_iCP_i, \qquad i=1,\ldots,n.
  \end{gather}

  As it was for the functors in Section~4, denote the restrictions
  $C|_{H_i}$, $i=1,\ldots,n$, by $C_i$. Then, as in Section~4, the
  operators $C_i$ map $H_i$ into $\tilde{H}_i$, that is,
  \begin{gather}\label{x26}
    C_i(H_i)\subset\tilde{H}_i,\qquad i=1,\ldots,n.
  \end{gather}
  It follows from~(\ref{x25}) and~(\ref{x26}) that
  \begin{gather}\label{x27}
    C\Gamma_i=\tilde{\Gamma}_iC_i,\qquad i=1,\ldots,n.
  \end{gather}
  The identities~(\ref{x25}) are equivalent to the inclusions
  $C(H_i)\subset\tilde{H}_i$, $i=1,\ldots,n$, whence it follows
  that
  \begin{gather}\label{x28}
    C_i=\tilde{\Gamma}_i^*C\Gamma_i,\qquad i=1,\ldots,n.
  \end{gather}

  Similarly to Section~4, identities~(\ref{x27}) allow to represent
  $C$ as
  \begin{gather}\label{x29}
    C=\frac1{\alpha}\sum_{i=1}^n\tilde{\Gamma}_iC_i\Gamma_i^*.
  \end{gather}

  The above presents all the similarities with the calculations performed
  in Section~4; the operator~$\hat{C}$ is now def\/ined dif\/ferently. For
  the operator $\hat{C}=\Diag(C_1,C_2,\ldots,C_n):\mathcal
  H\rightarrow\tilde{\mathcal H}$, it is easy to check that
  $\hat{C}Q_i=\tilde{Q}_i\hat{C}$, $i=1,\ldots,n$. Then
  $Q_i\hat{C}^*=\hat{C}^*\tilde{Q}_i$, $i=1,\ldots,n$. The latter
  allows to conclude that $\hat{C}^*(\Imp \tilde{Q}_i)\subset\Imp Q_i$
  and, consequently,
  \begin{gather}\label{x30}
    \hat{C}^*\tilde{Q}_i=Q_i\hat{C}^*\tilde{Q}_i,\qquad
    i=1,\ldots,n.
  \end{gather}

  Denote by $(\tilde{P}\hat{C}P)_{ij}$ the elements of the block
  matrix of the operator $\tilde{P}\hat{C}P:H_1\oplus H_2\oplus\cdots\oplus
H_n\rightarrow\tilde{H}_1\oplus\tilde{H}_2\oplus\cdots\oplus\tilde{H}_n$.
  Then $(\tilde{P}\hat{C}P)_{ij}=\frac1{\alpha^2}\sum\limits_{k=1}^n\tilde{\Gamma}_i^*
\tilde{\Gamma}_k C_k\Gamma_k^*\Gamma_j=
\frac1{\alpha^2}\tilde{\Gamma}_i^*
(\sum\limits_{k=1}^n\tilde{\Gamma}_k C_k\Gamma_k^*)\Gamma_j
=\frac1{\alpha}\tilde{\Gamma}_i^*C\Gamma_j=
\frac1{\alpha}\tilde{\Gamma}_i^*\tilde{\Gamma}_jC_j=(\tilde{P}\hat{C})_{ij}$,
  that is, $\tilde{P}\hat{C}P=\tilde{P}\hat{C}$ and, consequently,
  \begin{gather}\label{x31}
    \hat{C}^*\tilde{P}=P\hat{C}^*\tilde{P}.
  \end{gather}

  Identities~(\ref{x30}) and~(\ref{x31}) mean that
  $\hat{C}^*\in\Mor(\EuScript F'(\tilde{\pi}),\EuScript
  F'(\pi))$. Def\/ine $\EuScript F'(C)=\hat{C}^*$, and this f\/inishes the
  construction of the functor $\EuScript F'$.

\begin{lemma}\label{lem3}
  The functor $\EuScript F'$ is an equivalence between the categories.
\end{lemma}

\begin{proof}
  Let us show that the functor is univalent. Let
  $C,D\in\Mor(\pi,\tilde{\pi})$ and $C\neq D$, and show that $\EuScript
  F'(C)\neq \EuScript F'(D)$. Indeed, if $\EuScript F'(C)=
  \EuScript F'(D)$, i.e., $\hat{C}^*=\hat{D}^*$, then $C_i=D_i$, $\forall\,
  i=1,\ldots,n$. Let us use~(\ref{x29}),
  \[
  C=\frac1{\alpha}\sum_{i=1}^n\tilde{\Gamma}_iC_i\Gamma_i^*,\qquad
  D=\frac1{\alpha}\sum_{i=1}^n\tilde{\Gamma}_iD_i\Gamma_i^*.
  \]
  It follows from $C_i=D_i$, $i=1,\ldots,n$, and the form of the
  representation operators $C$ and $D$ that $C=D$ and, hence, the
  functor $\EuScript F'$ is univalent.

  Let us show that $\EuScript F'$ is complete. Let $R\in\Mor(\EuScript
  F'(\tilde{\pi}),\EuScript F'(\pi))$ and construct a linear operator
  in the set $\Mor(\pi,\tilde{\pi})$ such that the value of this
  functor on the morphism is $R$. It follows from $R\in\Mor(\EuScript
  F'(\tilde{\pi}),\EuScript F'(\pi))$ that the operator
  $R:\tilde{\mathcal H}\rightarrow\mathcal H$ satisf\/ies
  \[
  Q_iR\tilde{Q}_i=R\tilde{Q}_i, \qquad i=1,\ldots,n,\qquad
  PR\tilde{P}=R\tilde{P}.
  \]

  Denote by $\hat{r}$ an operator in $B(\mathcal H,\tilde{\mathcal
    H})$ such that $\hat{r}^*=R$. Then the latter identities can be
  rewritten as follows:
  \[
  Q_i\hat{r}^*\tilde{Q}_i=\hat{r}^*\tilde{Q}_i, \qquad i=1,\ldots,n,\qquad
  P\hat{r}^*\tilde{P}=\hat{r}^*\tilde{P},
  \]
  and, consequently,
  \begin{gather}\label{x32}
    \tilde{Q}_i\hat{r}Q_i=\tilde{Q}_i\hat{r},\qquad i=1,\ldots,n
  \end{gather}
  and
  \begin{gather}\label{x33}
    \tilde{P}\hat{r}P=\tilde{P}\hat{r}.
  \end{gather}

  Let now $r_{ij}$ be elements of the block matrix of the operator
  $\hat{r}$ from $H_1\oplus H_2\oplus\cdots\oplus H_n$ into
  $\tilde{H}_1\oplus\tilde{H}_2\oplus\cdots\oplus\tilde{H}_n$.
  Identities~(\ref{x32}) imply that if $i\neq j$, then $r_{ij}=0$.
  Denote $r_i=r_{ii}$, $i=1,\ldots,n$. Then
  $r_i:H_i\rightarrow\tilde{H}_i$, $i=1,\ldots,n$, and
  $\hat{r}=\Diag(r_1,r_2,\ldots,r_n)$. Consider
  $r:H\rightarrow\tilde{H}$ def\/ined by
  \begin{gather}\label{x34}
    r=\frac1{\alpha}\tilde{\Gamma}\hat{r}\Gamma^*.
  \end{gather}
  Identity~(\ref{x33}) and def\/inition~(\ref{x34}) imply that
  $\frac1{\alpha}\tilde{\Gamma}^*r\Gamma=\tilde{P}\hat{r}P=\tilde{P}\hat{r}$,
  then comparing the elements on the main diagonal of the
  corresponding block matrices gives
  \begin{gather}\label{x35}
    r_i=\tilde{\Gamma}_i^*r\Gamma_i,\qquad i=1,\ldots,n.
  \end{gather}

  Using the relation
  $\frac1{\alpha}\tilde{\Gamma}\tilde{\Gamma}^*=I_{\tilde{H}}$ we get
  $r\Gamma=I_{\tilde{H}}r\Gamma=
  (\frac1{\alpha}\tilde{\Gamma}\tilde{\Gamma}^*)r\Gamma=
  \tilde{\Gamma}(\frac1{\alpha}\tilde{\Gamma}^*r\Gamma)=
  \tilde{\Gamma}(\tilde{P}\hat{r}P)= \tilde{\Gamma}(\tilde{P}\hat{r})=
  \tilde{\Gamma}(\frac1{\alpha}\tilde{\Gamma}^*\tilde{\Gamma})\hat{r}=
  (\frac1{\alpha}\tilde{\Gamma}\tilde{\Gamma}^*)\tilde{\Gamma}\hat{r}=
  I_{\tilde{H}}\tilde{\Gamma}\hat{r}=\tilde{\Gamma}\hat{r}$.
  Rewrite the identity $r\Gamma=\tilde{\Gamma}\hat{r}$ in the matrix
  form,
  \[
  (r\Gamma_1\,r\Gamma_2\,\ldots r\Gamma_n)=
  (\tilde{\Gamma}_1r_1\,\tilde{\Gamma}_2r_2\,\ldots
  \,\tilde{\Gamma}_nr_n)
  \]
  that gives
  \begin{gather}\label{x36}
    r\Gamma_i=\tilde{\Gamma}_ir_i,\qquad i=1,\ldots,n.
  \end{gather}

  Using identities~(\ref{x35}) and (\ref{x36}) we get
  \begin{gather}\label{x37}
    rP_i=\tilde{P}_irP_i,\qquad i=1,\ldots,n.
  \end{gather}
  Indeed, $rP_i=r\Gamma_i\Gamma_i^*=\tilde{\Gamma}_ir_i\Gamma_i^*=
  \tilde{\Gamma}_i\tilde{\Gamma}_i^*r\Gamma_i\Gamma_i^*=\tilde{P}_irP_i$.
  Identities~(\ref{x37}) mean that $r\in\Mor(\pi,\tilde{\pi})$. Let us
  check that $\EuScript F'(r)=\hat{r}^*=R$. Denote by $C$ the
  constructed morphism $r$ and f\/ind $\EuScript F'(C)$. Since
  $\EuScript F'(C)=\hat{C}^*$, where
  $\hat{C}=\Diag(C_1,C_2,\ldots,C_n):\mathcal
  H\rightarrow\tilde{\mathcal H}$, let us f\/ind
  $C_i=C|_{H_i}=r|_{H_i}$, $i=1,\ldots,n$. Since
  $C\in\Mor(\pi,\tilde{\pi})$, it follows from~(\ref{x28})
  and~(\ref{x35}) that
  $C_i=\tilde{\Gamma}_i^*C\Gamma_i=\tilde{\Gamma}_i^*r\Gamma_i=r_i$.
  Then $\hat{C}=\hat{r}$, $\hat{C}^*=\hat{r}^*$ and $\EuScript
  F'(r)=\EuScript F'(C)=\hat{C}^*=\hat{r}^*=R$. This proves that the
  functor $\EuScript F'$ is complete.

  Since $\EuScript F'(\pi)=\EuScript F(\pi)$ for any
 $\pi\in\Rep\mathcal P_{n,\alpha}$ and $\EuScript F^2=\Id$, we see that
 $\EuScript F'$ is an equivalence of the categories.
\end{proof}

\begin{lemma}\label{lem4}
  If a system $S_{\pi}$, $\pi\in\Rep\mathcal P_{n,{\rm com}}$, of $n$
  subspaces is transitive, then the system $S_{\EuScript F(\pi)}$ of
  $n+1$ subspaces is transitive. Also, $S_{\pi}\cong S_{\tilde{\pi}}$
  if and only if $S_{\EuScript F(\pi)}\cong S_{\EuScript
    F(\tilde{\pi})}$.
\end{lemma}

\begin{proof}
  Since the functors $\EuScript F$ and $\EuScript F'$ coincide on the
  objects of the categories, the representations constructed using the
  functors and the corresponding systems of subspaces will coincide,
  $S_{\EuScript F(\pi)}=S_{\EuScript F'(\pi)}$ for $\forall\,
  \pi\in\Rep\mathcal P_{n,\alpha}$. Let $\pi,\tilde{\pi}\in
  \Rep\mathcal P_{n,\alpha}$, $\alpha\neq0$, and the systems of
  subspaces $S_{\pi}=(H;H_1,H_2,\ldots,H_n)$ and
  $S_{\tilde{\pi}}=(\tilde{H};\tilde{H}_1,\tilde{H}_2,\ldots,\tilde{H}_n)$,
  which correspond to the representations $\pi$ and $\tilde{\pi}$, be
  isomorphic, that is, $S_{\pi}\cong S_{\tilde{\pi}}$. By the
  def\/inition of isomorphic systems, there exists a linear operator $T\in
  B(H,\tilde{H})$ such that $T^{-1}\in B(\tilde{H},H)$ and
  $T(H_i)=\tilde{H}_i$, $i=1,\ldots,n$. It follows from
  $T(H_i)=\tilde{H}_i$, $i=1,\ldots,n$, that
  $T(H_i)\subset\tilde{H}_i$, $i=1,\ldots,n$, and, consequently,
  $TP_i=\tilde{P}_iTP_i$, $i=1,\ldots,n$. The latter identities
  mean that $T\in\Mor(\pi,\tilde{\pi})$. Then
  $\hat{T}^*\in\Mor(\EuScript F'(\tilde{\pi}),\EuScript F'(\pi))$ and
  \begin{gather}\label{x38}
    \hat{T}^*(\Imp \tilde{Q}_i)\subset(\Imp
    Q_i)\,(i=1,\ldots,n)\qquad \mbox{and} \qquad \hat{T}^*(\Imp
    \tilde{P})\subset(\Imp P).
  \end{gather}

  Again, considering the identities $T(H_i)=\tilde{H}_i$,
  $i=1,\ldots,n$, we conclude that $T(H_i)\supset\tilde{H}_i$,
  $i=1,\ldots,n$, that is, $T^{-1}(\tilde{H}_i)\subset H_i$,
  $i=1,\ldots,n$, and, respectively,
  $T^{-1}\tilde{P}_i=P_iT^{-1}\tilde{P}_i$, $i=1,\ldots,n$. These
  identities imply that $T^{-1}\in\Mor(\tilde{\pi},\pi)$. Then
  $\widehat{T^{-1}}^*\in\Mor(\EuScript F'(\pi),\EuScript
  F'(\tilde{\pi}))$, whence using
  $\widehat{T^{-1}}^*=(\hat{T}^{-1})^*=(\hat{T}^*)^{-1}$ we have
  \[
  \Imp \tilde{Q}_i\supset(\hat{T}^*)^{-1}(\Imp
  Q_i), \quad i=1,\ldots,n\qquad \mbox{and} \qquad \Imp
  \tilde{P}\supset(\hat{T}^*)^{-1}(\Imp P),
  \]
  and, consequently,
  \begin{gather}\label{x39}
    \hat{T}^*(\Imp \tilde{Q}_i)\supset \Imp Q_i,\quad
    i=1,\ldots,n\qquad \mbox{and} \qquad \hat{T}^*(\Imp
    \tilde{P})\supset\Imp P.
  \end{gather}
  It follows from~(\ref{x38}) and~(\ref{x39}) that
  \[
  \hat{T}^*(\Imp
  \tilde{Q}_i)= \Imp Q_i, \quad i=1,\ldots,n\qquad \mbox{and} \qquad
  \hat{T}^*(\Imp \tilde{P})=\Imp P
  \]
  that shows that it is an isomorphism of the systems that
  correspond to the representations~$\EuScript F'(\pi)$ and $\EuScript
  F'(\tilde{\pi})$ and, since the functors $\EuScript F'$ and
  $\EuScript F$ coincide on the objects of the category, it is an
  isomorphism of the systems corresponding to the representations $\EuScript
  F(\pi)$ and $\EuScript F(\tilde{\pi})$.

  Since the functor $\EuScript F'$ is complete, similar reasonings
  show that the representations that correspond to nonisomorphic
  systems are mapped by the functor $\EuScript F'$, and hence the
  functor $\EuScript F$, into representations that give rise to
  nonisomorphic systems.

  Let us now prove the f\/irst part of the proposition. Let $\pi$ be a
  representation of the algebra~$\mathcal P_{n,\alpha}$ and
  $\pi(p_i)=P_i$, $i=1,\ldots,n$, be orthogonal projections on a
  representation space $H$. And let the system of orthogonal
  projections $P_1, P_2, \ldots, P_n$ induce a transitive
  system of subspaces $S_{\pi}=(H;H_1,H_2,\ldots,H_n)$, where
  $H_i=P_iH$, $i=1,\ldots,n$, that is,
  \[
  \End(S_{\pi})=\{r\in B(H)\,|\,r(H_i)\subset H_i,
  i=1,\ldots,n\}=\Mor(\pi,\pi)=\mathbb CI.
  \]

  Consider $\EuScript F'(\pi)=\hat{\pi}$, where $\hat{\pi}(q_i)=Q_i$,
  $i=1,\ldots,n$, and $\hat{\pi}(p)=P$, and the corresponding
  system $S_{\hat{\pi}}$ of subspaces. Let now $R\in
  \End(S_{\hat{\pi}})$. Using $\End(S_{\hat{\pi}})=\Mor(\EuScript
  F'(\pi),\EuScript F'(\pi))$ and since the functor $\EuScript F'$ is
  complete, we see that $\EuScript F'(r)=R$, where $r\in\Mor(\pi,\pi)$
  is constructed from the diagonal operator
  $R^*=\Diag(r_1,r_2,\ldots,r_n)$ on the space $H_1\oplus
  H_2\oplus\cdots\oplus H_n$ as follows:
  \begin{gather}\label{x40}
    r=\frac1{\alpha}\Gamma R^*\Gamma^*.
  \end{gather}

  Using the inclusion $R\in \Mor(\EuScript F'(\pi),\EuScript
  F'(\pi))$, which is similar to identity~(\ref{x35}), we get
  \begin{gather}\label{x41}
    r_i=\Gamma_i^*r\Gamma_i,\qquad i=1,\ldots,n.
  \end{gather}

  Since the system $S_{\pi}$ is transitive, the operator $r$ is
  scalar, that is, $r=\lambda I_{H}$. Using
  $\Gamma_i^*\Gamma_i=I_{H_i}$, $i=1,\ldots,n$, and
  identities~(\ref{x41}) we get
  \[
  r_i=\lambda I_{H_i},\qquad i=1,\ldots,n.
  \]
  Then $R^*$ is a scalar operator and, consequently, $R$ is also a
  scalar operator that means that the system $S_{\EuScript F'(\pi)}$
  is transitive and such is $S_{\EuScript F(\pi)}$.
\end{proof}

The claim of Theorem~\ref{t4} follows from Theorem~\ref{t3} and
Lemma~\ref{lem4}.

\subsection{Transitive quintuples of subspaces}

By Theorem~\ref{t4}, the functor $\EuScript F$ maps known
nonisomorphic transitive quadruples of subspaces of the form
$S_{\pi}$, where $\pi\in\Rep\mathcal P_{4,{\rm com}}$, into
nonisomorphic transitive quintuples $S_{\EuScript
F(\pi)}$~\cite{MS,M}. In this section, we give inequivalent
irreducible $*$-representations of the $*$-algebras $\mathcal
P_{4,{\rm abo},\tau}$, $\tau\in\tilde{\Sigma}_4$, where
$\tilde{\Sigma}_4$ that is the set of $\tau\in\mathbb R_+$ such
that there exists at least one $*$-representation of the
$*$-algebra $\mathcal P_{4,{\rm abo},\tau}$, is related to $\Sigma_4$,
the set $\alpha\in\mathbb R_+$ such that there exists at least one
$*$-representation of the $*$-algebra $\mathcal P_{4,\alpha}$, via
the following relation~\cite{PS}:
\[
\tilde{\Sigma}_4=\{0\}\cup\left\{\frac1{\alpha}\,
|\,\alpha\neq0,\alpha\in\Sigma_4\right\}.
\]
Here, by~\cite{KRS}, $\Sigma_4=\big\{0,1,2-\frac2{2k+1}\,
(k=1,2,\ldots), 2-\frac1{n}\, (n=2,3,\ldots), 2,2+\frac1{n}\,
(n=2,3,\ldots),2+\frac2{2k+1}\, (k=1,2,\ldots),3,4\big\} $. For
these representations, the corresponding systems of subspaces are
nonisomorphic and transitive.

Let $e_{i,j}^{r\times s}$ denote an ($r\times s$)-matrix that has
$1$ at the intersection of the $i$th row and the $j$th column,
with other elements being zero.

1) The $*$-algebra $\mathcal P_{4,{\rm abo},0}$ has $4$ irreducible
inequivalent one-dimensional representations, $Q_1=\cdots=
Q_{k-1}=Q_{k+1}=\cdots=Q_4=P=0$, $Q_k=1$.

2) The $*$-algebra $\mathcal P_{4,{\rm abo},1}$ has $4$ irreducible
inequivalent one-dimensional representations, $Q_1=\cdots=
Q_{k-1}=Q_{k+1}=\cdots=Q_4=0$, $P=Q_k=1$.

3) The $*$-algebra $\mathcal P_{4,{\rm abo},\frac13}$ has $4$
irreducible inequivalent three-dimensional representations that
are unitary equivalent, up to a permutation, to
\begin{gather*}
Q_1=1\oplus0\oplus0,\qquad
 Q_3=0\oplus0\oplus1,\qquad  P=\frac13\sum_{i,j=1}^3e_{i,j}^{3\times 3},\\
 Q_2=0\oplus1\oplus0,\qquad Q_4=0\oplus0\oplus0,\qquad \mathcal H=\mathbb C\oplus\mathbb
C\oplus\mathbb C.
 \end{gather*}

4) The $*$-algebra $\mathcal P_{4,{\rm abo},\frac14}$ has a unique
irreducible four dimensional representation,
\begin{gather*}
Q_1=1\oplus0\oplus0\oplus0,\qquad
Q_3=0\oplus0\oplus1\oplus0,\qquad P=\frac14\sum_{i,j=1}^4e_{i,j}^{4\times 4},\\
 Q_2=0\oplus1\oplus0\oplus0 , \qquad
 Q_4=0\oplus0\oplus0\oplus1,\qquad \mathcal H=\mathbb C\oplus\mathbb
C\oplus\mathbb C\oplus\mathbb C.
 \end{gather*}

5) The $*$-algebra $\mathcal P_{4,{\rm abo},\frac12}$ has $6$
irreducible two-dimensional representations that are unitary
equivalent, up to a permutation, to
\[
Q_1=1\oplus0,\qquad Q_2=0\oplus1,\qquad Q_3=0\oplus0,\qquad
Q_4=0\oplus0,\qquad P=\frac12\begin{pmatrix}1&1\\1&1\end{pmatrix},
\]
where the representation space is $\mathcal H=\mathbb
C\oplus\mathbb C$, and the following inequivalent four
dimensional representations that depend on the points of the set
$\Omega=\{(a,b,c)\in\mathbb R| a^2+b^2+c^2=1, a>0$, $b>0,
c\in(-1,1);\; \mbox{or}\; a=0, b^2+c^2=1, b>0, c>0 ;\; \mbox{or}\;
b=0, a^2+c^2=1, a>0, c>0\}$:
\begin{gather*}
Q_1=1\oplus0\oplus0\oplus0,\qquad
 Q_3=0\oplus0\oplus1\oplus0,\\
 Q_2=0\oplus1\oplus0\oplus0,\qquad
 Q_4=0\oplus0\oplus0\oplus1,
\\
P=\frac12
\begin{pmatrix}1&\frac{c(c-ib)}{\sqrt{1-a^2}}&
  \frac{b(b+ic)}{\sqrt{1-a^2}}
  &a\\[2mm]
  \frac{c(c+ib)}{\sqrt{1-a^2}}&1&-a&\frac{b(b-ic)}{\sqrt{1-a^2}}
  \\[2mm]
  \frac{b(b-ic)}{\sqrt{1-a^2}}&-a&1&\frac{c(c+ib)}{\sqrt{1-a^2}}
  \\[2mm]
  a&\frac{b(b+ic)}{\sqrt{1-a^2}}
  &\frac{c(c-ib)}{\sqrt{1-a^2}}&1\end{pmatrix},
\end{gather*}
where the representation space is $\mathcal H=\mathbb
C\oplus\mathbb C\oplus\mathbb C\oplus\mathbb C$.

6) The $*$-algebras $\mathcal P_{4,{\rm abo},\frac1{\alpha}}$, for
$\alpha=2-\frac2{2k+1}$, $k=1,2,\ldots$, have unique irreducible
representations
\begin{gather*}
Q_1=I\oplus0\oplus0\oplus0,\qquad
Q_3=0\oplus0\oplus I\oplus0,\\
 Q_2=0\oplus I\oplus0\oplus0 ,\qquad
 Q_4=0\oplus0\oplus0\oplus I,
 \\
P=\frac1{\alpha}\begin{pmatrix}A&B\\B^t&C\end{pmatrix},\qquad
\mbox{where}\quad A=\begin{pmatrix}I&A_1\\A_1&I\end{pmatrix},\quad
C=\begin{pmatrix}I&C_1\\C_1&I\end{pmatrix},\quad
B=\begin{pmatrix}B_{00}&B_{01}\\B_{10}&B_{11}\end{pmatrix},
\\
A_1=\frac1{2k+1}\sum_{i=1}^k(2k+3-4i)e_{i,i}^{k\times k},\qquad
C_1=\frac1{2k+1}\sum_{i=1}^k(2k+1-4i)e_{i,i}^{k\times k},
\\
B_{lm}=\frac{(-1)^\ell}{2k+1}\sum_{i=1}^k\sqrt{(2k-2i+1)(2i-1)}\,e_{i,i}^{k\times
k}+\frac{(-1)^m}{2k+1}\sum_{i=1}^{k-1}\sqrt{(2k-2i)2i}\,e_{i+1,i}^{k\times
k},
\end{gather*}
and the representation space is
\[
\mathcal H= \mathbb C^{k}\oplus\mathbb C^{k}\oplus\mathbb
C^{k}\oplus\mathbb C^{k}.
\]

7) The $*$-algebras $\mathcal P_{4,{\rm abo},\frac1{\alpha}}$, for
$\alpha=2-\frac1{2k+1}$, $k=1,2,\ldots$, have  unique irreducible
representations
\begin{gather*}
Q_1=I\oplus0\oplus0\oplus0,\qquad
Q_3=0\oplus0\oplus I\oplus0,\\
Q_2=0\oplus I\oplus0\oplus0 ,\qquad Q_4=0\oplus0\oplus0\oplus I,
 \\
P=\frac1{\alpha}\begin{pmatrix}A&B\\B^t&C\end{pmatrix},\qquad
\mbox{where}\quad
A=\begin{pmatrix}I&A_1\\A_1^t&I\end{pmatrix},\quad
C=\begin{pmatrix}I&C_1\\C_1&I\end{pmatrix},\quad
B=\begin{pmatrix}\eta&\eta\\B_{00}&B_{10}\\B_{01}&B_{11}\end{pmatrix},
\\
\eta=({\scriptstyle\sqrt{\frac{k}{2k+1}}},\underbrace{0,0,\ldots,0}_{k-1}),
\\
A_1=-\frac1{2k+1}\sum_{i=1}^k2ie_{i+1,i}^{(k+1)\times k},\qquad
C_1=-\frac1{2k+1}\sum_{i=1}^k(2i-1)e_{i,i}^{k\times k},
\\
B_{lm}=\frac{(-1)^\ell}{4k+2}\sum_{i=1}^k
\!\sqrt{(2k-2i+1)(2k+2i)}\,e_{i,i}^{k\times k}\!
+\frac{(-1)^m}{4k+2}\sum_{i=1}^{k-1}\!\sqrt{(2k-2i)(2k+2i+1)}\,e_{i,i+1}^{k\times
k},\!
\end{gather*}
and the representation space is
\[
\mathcal H= \mathbb C^{k+1}\oplus\mathbb C^{k}\oplus\mathbb
C^{k}\oplus\mathbb C^{k}.
\]

8) The $*$-algebras $\mathcal P_{4,{\rm abo},\frac1{\alpha}}$, for
$\alpha=2-\frac1{2k}$, $k=1,2,\ldots$, have  unique irreducible
representations
\begin{gather*}
Q_1=I\oplus0\oplus0\oplus0,\qquad
Q_3=0\oplus0\oplus I\oplus0,\\
Q_2=0\oplus I\oplus0\oplus0 ,\qquad Q_4=0\oplus0\oplus0\oplus I,
 \\
P=\frac1{\alpha}\begin{pmatrix}A&B\\B^t&C\end{pmatrix},\qquad
\mbox{where}\quad
A=\begin{pmatrix}I&A_1\\A_1^t&I\end{pmatrix},\quad
C=\begin{pmatrix}I&C_1\\C_1&I\end{pmatrix},\quad
B=\begin{pmatrix}B_{00}&B_{10}\\\eta&\eta\\B_{01}&B_{11}\end{pmatrix},
\\
\eta=({\scriptstyle\sqrt{\frac{2k-1}{4k}}},\underbrace{0,0,\ldots,0}_{k-1}),
\\
A_1=-\frac1{k}\sum_{i=1}^{k-1}ie_{i,i+1}^{(k-1)\times k},\qquad
C_1=-\frac1{2k}\sum_{i=1}^k(2i-1)e_{i,i}^{k\times k},
\\
B_{lm}=\frac{(-1)^\ell}{4k}\sum_{i=1}^{k-1}\sqrt{(2k-2i)(2k+2i-1)}\,e_{i,i}^{(k-1)\times
k}\\
\phantom{B_{lm}=}{}
+\frac{(-1)^m}{4k}\sum_{i=1}^{k-1}\sqrt{(2k-2i-1)(2k+2i)}\,e_{i,i+1}^{(k-1)\times
k},
\end{gather*}
and the representation space is
\[
\mathcal H= \mathbb C^{k-1}\oplus\mathbb C^{k}\oplus\mathbb
C^{k}\oplus\mathbb C^{k}.
\]

9) The $*$-algebras $\mathcal P_{4,{\rm abo},\frac1{\alpha}}$, for
$\alpha=2+\frac1{2k}$, $k=1,2,\ldots$, have  unique irreducible
representations
\begin{gather*}
Q_1=I\oplus0\oplus0\oplus0,\qquad
Q_3=0\oplus0\oplus I\oplus0,\\
Q_2=0\oplus I\oplus0\oplus0 ,\qquad Q_4=0\oplus0\oplus0\oplus I,
\\
P=\frac1{\alpha}\begin{pmatrix}A&B\\B^t&C\end{pmatrix},\qquad
\mbox{where}\quad
A=\begin{pmatrix}I&A_1\\A_1^t&I\end{pmatrix},\quad
C=\begin{pmatrix}I&C_1\\C_1&I\end{pmatrix},\quad
B=\begin{pmatrix}\eta&\eta\\B_{11}&B_{01}\\B_{10}&B_{00}\end{pmatrix},
\\
\eta=({\scriptstyle\sqrt{\frac{2k+1}{4k}}},\underbrace{0,0,\ldots,0}_{k-1}),
\\
A_1=\frac1{k}\sum_{i=1}^{k}ie_{i+1,i}^{(k+1)\times k},\qquad
C_1=\frac1{2k}\sum_{i=1}^k(2i-1)e_{i,i}^{k\times k},
\\
B_{lm}=\frac{(-1)^\ell}{4k}\sum_{i=1}^{k}\sqrt{(2k+2i)(2k-2i+1)}\,e_{i,i}^{k\times
k}\\
\phantom{B_{lm}=}{}
+\frac{(-1)^m}{4k}\sum_{i=1}^{k-1}\sqrt{(2k+2i+1)(2k-2i)}\,e_{i,i+1}^{k\times
k},
\end{gather*}
and the representation space is
\[
\mathcal H= \mathbb C^{k+1}\oplus\mathbb C^{k}\oplus\mathbb
C^{k}\oplus\mathbb C^{k}.
\]

10) The $*$-algebras $\mathcal P_{4,{\rm abo},\frac1{\alpha}}$, for
$\alpha=2+\frac1{2k+1}$, $k=1,2,\ldots$, have unique irreducible
representations
\begin{gather*}
Q_1=I\oplus0\oplus0\oplus0,\qquad
Q_3=0\oplus0\oplus I\oplus0,\\
Q_2=0\oplus I\oplus0\oplus0 ,\qquad Q_4=0\oplus0\oplus0\oplus I,
 \\
P=\frac1{\alpha}\begin{pmatrix}A&B\\B^t&C\end{pmatrix},\qquad
\mbox{where}\quad
A=\begin{pmatrix}I&A_1\\A_1^t&I\end{pmatrix},\quad
C=\begin{pmatrix}I&C_1\\C_1&I\end{pmatrix},\quad
B=\begin{pmatrix}B_{11}&B_{01}\\\eta&\eta\\B_{10}&B_{00}\end{pmatrix},
\\
\eta=({\scriptstyle\sqrt{\frac{k+1}{2k+1}}},\underbrace{0,0,\ldots,0}_{k}),
\\
A_1=\frac1{2k+1}\sum_{i=1}^{k}2ie_{i,i+1}^{k\times (k+1)},\qquad
C_1=\frac1{2k+1}\sum_{i=1}^{k+1}(2i-1)e_{i,i}^{(k+1)\times (k+1)},
\\
B_{lm}=\frac{(-1)^\ell}{4k+2}\sum_{i=1}^{k}\sqrt{(2k-2i+2)(2k+2i+1)}\,e_{i,i}^{k\times
(k+1)}\\
\phantom{B_{lm}=}{}
+\frac{(-1)^m}{4k+2}\sum_{i=1}^{k}\sqrt{(2k+2i-1)(2k+2i+2)}\,e_{i,i+1}^{k\times
(k+1)},
\end{gather*}
and the representation space is
\[
\mathcal H= \mathbb C^{k}\oplus\mathbb C^{k+1}\oplus\mathbb
C^{k+1}\oplus\mathbb C^{k+1}.
\]

11) The $*$-algebras $\mathcal P_{4,{\rm abo},\frac1{\alpha}}$, for
$\alpha=2+\frac2{2k+1}$, $k=1,2,\ldots$, have unique irreducible
representations
\begin{gather*}
Q_1=I\oplus0\oplus0\oplus0,\qquad
Q_3=0\oplus0\oplus I\oplus0,\\
Q_2=0\oplus I\oplus0\oplus0 ,\qquad Q_4=0\oplus0\oplus0\oplus I,
\\
P=\frac1{\alpha}\begin{pmatrix}A&B\\B^t&C\end{pmatrix},\qquad
\mbox{where}\quad A=\begin{pmatrix}I&A_1\\A_1&I\end{pmatrix},\quad
C=\begin{pmatrix}I&C_1\\C_1&I\end{pmatrix},\quad
B=\begin{pmatrix}B_{11}&B_{01}\\B_{10}&B_{00}\end{pmatrix},
\\
A_1=-\frac1{2k+1}\sum_{i=1}^{k+1}(2k+3-4i)e_{i,i}^{(k+1)\times
(k+1)},
\\
C_1=e_{1,1}^{(k+1)\times(k+1)}-\frac1{2k+1}\sum_{i=2}^{k+1}(2k+5-4i)e_{i,i}^{(k+1)\times
(k+1)},
\\
B_{lm}=\frac1{\sqrt{2k+1}}e_{1,1}^{(k+1)\times(k+1)}
+\frac{(-1)^\ell}{2k+1}\sum_{i=2}^{k+1}\sqrt{(2k-2i+3)(2i-1)}\,e_{i,i}^{(k+1)\times
(k+1)}\\
\phantom{B_{lm}=}{}
+\frac{(-1)^m}{2k+1}\sum_{i=1}^{k}\sqrt{(2k-2i+2)2i}\,e_{i,i+1}^{(k+1)\times
(k+1)},
\end{gather*}
and the representation space is
\[
\mathcal H= \mathbb C^{k+1}\oplus\mathbb C^{k+1}\oplus\mathbb
C^{k+1}\oplus\mathbb C^{k+1}.
\]

\subsection*{Acknowledgments}
The authors are grateful to S.A.~Kruglyak for useful remarks and
suggestions.

\LastPageEnding


\begin{thebibliography}{99}
\footnotesize


\bibitem{Ha1} Halmos P.R., Two subspaces, {\it Trans. Amer. Math.
    Soc.}, 1969, V.144, 381--389.

\bibitem{Ha2} Halmos P.R., Ten problems in Hilbert space, {\it Bull.
    Amer. Math. Soc.}, 1970, V.76, 887--933.

\bibitem{B} Brenner S., Endomorphism algebras of vector spaces with
  distinguished sets of subspaces, {\it J. Algebra}, 1967, V.6,
  100--114.

\bibitem{GP} Gel'fand I.M., Ponomarev V.A., Problems of linear algebra
  and classif\/ication of quadruples of subspaces in a
  f\/inite-dimensional vector space, {\it Coll. Math. Spc. Bolyai 5, Tihany},
  1970, 163--237.

\bibitem{N} Nazarova L.A., Representations of a quadruple, {\it Izv.
    AN SSSR}, 1967, V.31, N~6, 1361--1377 (in Russian).

\bibitem{EW} Enomoto M., Watatani Ya., Relative position of four
  subspaces in a Hilbert space, math.OA/0404545.

\bibitem{KRS} Kruglyak S.A., Rabanovich V.I., Samo\v\i{}lenko
  Yu.S., On sums of projections, {\it Funktsional. Anal. i Prilozhen.}, 2002, V.36, N~3, 30--35  (English transl.:
  {\it Funct. Anal. Appl.}, 2002, V.36, N~3, 182--195).

\bibitem{KS} Kruglyak S.A., Samo\v\i{}lenko Yu.S., On the complexity of
  description of representations of $*$-algebras generated by
  idempotents, {\it Proc. Amer. Math. Soc.}, 2000, V.128, 1655--1664.

\bibitem{OS} Ostrovskyi V.L., Samo\v\i{}lenko Yu.S., Introduction to the theory
  of representations of f\/initely presented $*$-al\-gebras. I.
  Representations by bounded operators, Harwood Acad. Publs., 1999.

\bibitem{PS} Popovich S.V., Samo\v\i{}lenko Yu.S., On homomorphisms of
  algebras generated by projections, and the Coxeter functors,  {\it
    Ukrain. Mat. Zh.}, 2003, V.55, N~9, 1224--1237 (English transl.: {\it Ukrainian Math. J.}, 2003,
    V.55, N~9, 1480--1496).

\bibitem{MS} Moskaleva Yu.P., Samo\v\i{}lenko Yu.S., Systems of $n$
  subspaces and representations of $*$-algebras generated by
  projections, {\it Methods Funct. Anal. Topology}, 2006, V.12, N~1, 57--73.

\bibitem{M} Moskaleva Yu.P., On $*$-representations of the algebra
  $\mathcal P_{4,{\rm abo},\tau}$, {\it Uchenye Zapiski Tavricheskogo Natsional'nogo Universiteta
  imeni Vernadskogo, Seriya Matem. Mech. Inform. Kibern.}, 2005,
  N~1, 27--35 (in Russian).

\end{thebibliography}
\end{document}